\documentclass[11pt]{article}
\usepackage{amsmath,amssymb,euscript,graphicx,amsfonts}
\usepackage{color}
\usepackage[dvips]{epsfig}
\usepackage{hyperref}

\title{HAMILTONICITY OF VERTEX-TRANSITIVE GRAPHS OF ORDER $6p$}
\date{\today}

\newtheorem{theorem}{Theorem}[section]
\newtheorem{proposition}[theorem]{Proposition}
\newtheorem{lemma}[theorem]{Lemma}

\newcommand{\Qed}{\rule{2.5mm}{3mm}}

\newcommand{\Aut}{\hbox{{\rm Aut}}}

\newcommand{\p}{\wp}

\newcommand{\la}{\langle}
\newcommand{\ra}{\rangle}

\newcommand{\ZZ}{\mathbb{Z}}

\newcommand{\AGL}{\hbox{{\rm AGL}}}

\renewcommand{\S}{{\mathcal S}}
\renewcommand{\P}{{\cal P}}

\renewcommand{\O}{{\mathcal{O}}}

\newcommand{\B}{{\cal{B}}}

\newcommand{\R}{\mathcal{R}}
\newcommand{\M}{\mathcal{M}}

\newcommand{\PSL}{\hbox{{\rm PSL}}}

\newcommand{\proofT}{ {\sc Proof of Theorem~\ref{the:main}: }}

\newenvironment{proof}{{\noindent \sc Proof.}}{\hfill $\Qed$ \\}

 \newcounter{case}
 \newenvironment{case}[1][\unskip]{\refstepcounter{case}\sc
 \smallskip \noindent Case \thecase\ #1.\ }{\unskip\upshape}
 \renewcommand{\thecase}{\arabic{case}}
 
 \newcounter{subcase}
 
\numberwithin{subcase}{case}

\def\ZZ{{\hbox{\sf Z\kern-.43emZ}}}

\input amssym.def

\input amssym.tex

\color{black}
\begin{document}


\begin{center}
{\bf\large HAMILTON PATHS AND CYCLES IN VERTEX-TRANSITIVE GRAPHS OF ORDER $6p$}
\end{center}
\bigskip\noindent
\begin{center}

  Klavdija Kutnar{\small$^{a,}$}\footnotemark \, 
 and Primo\v z \v Sparl{\small$^b$}$^,$*\\

\bigskip

{\it {\small  $^a$University of Primorska, Cankarjeva 6, 6000 Koper, Slovenia\\
$^b$IMFM, University of Ljubljana, Jadranska 19, 1000 Ljubljana, Slovenia}}
\
\end{center}

\addtocounter{footnote}{0}  \footnotetext{Supported in part by
``Agencija za raziskovalno dejavnost Republike Slovenije'', research program P1-0285.

~*Corresponding author e-mail: ~primoz.sparl@fmf.uni-lj.si}

\begin{abstract}
 It is shown that every connected vertex-transitive graph of order $6p$,
where $p$ is a prime, contains a Hamilton path. Moreover, it is shown that, except for the truncation of the
Petersen graph, every connected vertex-transitive graph of order $6p$ which is not genuinely imprimitive 
contains a Hamilton cycle.
\end{abstract}

\bigskip
\begin{quotation}
\noindent {\em Keywords:} graph, vertex-transitive, Hamilton
cycle,  Hamilton path, automorphism group.
\end{quotation}


\bigskip
\section{Introductory remarks}
\label{sec:intro}
\indent

This  paper deals with the existence of Hamilton paths and Hamilton cycles in connected vertex-transitive
graphs of order $6p$, where $p$ is a prime.
(Throughout this paper $p$ will always denote a prime number.)
The question
whether every connected vertex-transitive graph contains a Hamilton path was posed by 
Lov\'asz in 1969 (see \cite{LL70}). So far no example giving a negative answer to this question 
has been found.
Moreover, apart from the trivial example $K_2$, there are only four known connected vertex-transitive graphs, 
which do not contain a  Hamilton cycle. 
These are the Petersen graph, the Coxeter graph and the truncations of these two graphs, 
that is the graphs obtained from them by replacing each vertex by a triangle. 
This supports the conjecture of Thomassen~\cite{BW78,CT91} that only finitely many connected
vertex-transitive graphs without a Hamilton cycle exist. On the other hand, Babai~\cite{LB79, HC95} 
conjectured that infinitely many such graphs exist.

Despite the fact that these questions have been challenging mathematicians for almost forty
years, only partial results have been obtained thus far. 
For instance, it is known that connected vertex-transitive graphs of orders
$kp$, where $k \leq 5$, $p^j$, where $j \leq 4$, and $2p^2$ contain a  Hamilton path.
Furthermore, for all of these families, except for the graphs of order $5p$, it is
also known that they contain a Hamilton cycle (except for the above mentioned
Petersen and Coxeter graph), see \cite{BA79,YC98, KM??, DM85, DM87,DM88,MP82,MP83, T67}.
The problem has also been considered for the subclass of Cayley graphs, resulting in a 
number of partial results (see for example ~\cite{AZ89,CG96,GM06,KW85,DM83, DW85,DW86}).
Also, it is known that every connected vertex-transitive graph, other than the Petersen graph,
whose automorphism group contains a transitive subgroup with a cyclic  
commutator subgroup of prime-power order, has a Hamilton cycle.
The result was proved in~\cite{DGMW98} and it uses a results from a series of papers
dealing with the same group-theoretic restrictions in the context of Cayley graphs
\cite{ED83,DM83,DW85}.

The main object of this paper is to show that every
connected vertex-transitive graph of order $6p$ contains a
Hamilton path. This result represents a new building block of the project to show that
all connected vertex-transitive graphs on up to $100$ vertices have this property.

\begin{theorem}
\label{the:main}
	Every connected vertex-transitive graph of order $6p$, where $p$ is a prime, contains a Hamilton path. 
	Moreover, with the exception of the truncation of the Petersen graph, every such
	graph which is not genuinely imprimitive 
	contains a Hamilton cycle.
\end{theorem}

The paper is organized as follows. In Section~\ref{sec:notation} notions concerning this paper
are introduced together with the notation and some auxiliary results that are needed
in the subsequent sections. The rest of the paper is devoted to proving Theorem~\ref{the:main}.
As a vertex-transitive graph is either 
genuinely imprimitive, quasiprimitive
or primitive, we divide our investigation depending on which of these three families the graph in question
belongs to.
The genuinely imprimitive graphs are considered in Section~\ref{sec:three}. The 
investigation of these graphs depends on the size of the corresponding blocks.
As for the quasiprimitive and primitive graphs of order $6p$, they are known (see \cite{GP, MSZ95}).
Therefore, the existence of Hamilton paths (or cycles) in these graphs can (at least in general) be verified. 
This is done in Sections~\ref{sec:quasi} and \ref{sec:pri}.
Finally, the results are combined in Section~\ref{sec:proof}, where the Theorem~\ref{the:main}
is proved.


\section{Notation and preliminary results}
\label{sec:notation}
\indent

Throughout this paper graphs are finite, simple and undirected, and groups are finite, unless specified otherwise. 
Furthermore, a {\em multigraph} is a generalization of a graph in which we allow multiedges and loops.
Given a graph $X$ we let $V(X)$ and $E(X)$ be the vertex set and the edge set of $X$, respectively.
For adjacent vertices $u,v \in V(X)$ we write $u \sim v$ and denote the corresponding edge by $uv$.
Let $U$ and $W$ be disjoint subsets of $V(X)$.
The subgraph of $X$ induced by $U$ will be denoted by $X\la U \ra$.
Similarly, we let $[U,W]$  denote the bipartite subgraph
of $X$ induced by the edges having one endvertex in $U$
and the other endvertex in $W$.

Given a transitive group $G$ acting on a set $V$,
we say that a partition $\B$ of $V$ is $G$-{\em invariant}
if the elements of $G$ permute the parts, that is, {\em blocks} of $\B$, setwise.
If the trivial partitions $\{V\}$ and $\{\{v\}: v \in V\}$ are the only
$G$-invariant partitions of $V$, then $G$ is said to be {\em primitive},
and is said to be {\em imprimitive} otherwise.
In the latter case we shall refer to a corresponding
$G$-invariant partition as to a  {\em complete imprimitivity block system}, in short 
an {\em imprimitivity block system}, of $G$.

A graph $X$ is said to be {\em vertex-transitive}
if its automorphism group, denoted by $\Aut X$, acts transitively on $V(X)$.
A vertex-transitive graph for which each transitive subgroup of its automorphism group is primitive
is called a {\em primitive graph}. Otherwise it is called an {\em imprimitive graph}.
If $X$ is imprimitive with an imprimitivity block system which is
formed by the orbits of a normal subgroup of some transitive subgroup $G \leq \Aut X$,
then the graph $X$ is said to be {\em genuinely imprimitive}.
If $X$ is imprimitive, but there exists no transitive subgroup $G$ of the 
automorphism group of $X$ having a nontransitive normal subgroup, then $X$ 
is said to be {\em quasiprimitive}.
Note that if $\B$ is an imprimitivity block system of some vertex-transitive graph, then
any two blocks $B,B' \in \B$ induce isomorphic vertex-transitive
subgraphs. 

The following simple observation about imprimitive groups of certain degrees will be useful latter on.

\begin{lemma}
\label{lem:blocksemi}
Let $G$ be an imprimitive permutation group of degree $mq$, $q$ a prime, with
a complete imprimitivity block system $\B$ and let $H \leq G$ have $m$ orbits 
of length $q$. Let $S$ be an orbit of $H$ and let $B \in \B$ be such that $B \cap S \neq \emptyset$. 
Then one of the following holds:
\begin{itemize}
	\item[(i)] $|B \cap S| = 1$, in which case $|B \cap S'| = 1$ for every 
	orbit $S'$ of $H$ which meets $B$, or 
	\item[(ii)] $B \cap S = S$, in which case $q$ divides $|B|$.
\end{itemize}
\end{lemma}

\begin{proof}
Let us first show that $|B \cap S|$ equals either to $1$ or to $q$. Suppose there exist distinct points
$u,v \in B \cap S$.
As $S$ is of prime length $q$, there exists some $\varphi \in H$, mapping $u$ to $v$, such that the restriction of $\varphi$
to $S$, denoted by $\varphi|_S$, is of order $q$. Then 
the orbit of $\varphi$ containing $u$ coincides with $S$. As $u\varphi = v$ and $u,v \in B$,
the block $B$ is fixed by $\varphi$. Consequently, $S \subseteq B$.

Suppose now that $B \cap S = \{u\}$ but $B \cap S' = S'$ for some orbits $S$ and $S'$ of $H$. 
In view of $B \cap S \neq S$, some element of $H$ moves the block $B$ to some other block. On the other hand (as
$B \cap S' = S'$), every element of $H$ fixes $B$ setwise. This contradiction proves $(i)$.
That $q$ divides $|B|$ when $B \cap S = S$ is now clear.
\end{proof}

Given a graph $X$ and a partition $\P$ of its vertex set we let the {\em quotient graph corresponding to
$\P$} be the graph $X_\P$ whose vertex set equals $\P$ with $A, B \in \P$ adjacent if there
exist vertices $a \in A$ and $b \in B$, such that $a \sim b$ in $X$. 

Let $m\geq 1$ and $n\geq 2$ be integers.
An automorphism of a graph is called $(m,n)$-{\em semiregular} if it has 
$m$ orbits of length $n$ and no other orbit.
Let now $X$ be a graph admitting an $(m,n)$-semiregular automorphism $\rho$ and denote the set of the orbits of $\rho$
by $\S$. Let $S, S' \in \S$. Clearly, the graph $[S,S']$ is regular. We let $d(S,S')$ denote
the valency of $[S,S']$. We let
the {\em quotient multigraph corresponding to $\rho$} be the multigraph $X_\rho$ whose vertex set
is $\S$ and in which $S, S' \in \S$ are joined by $d(S,S')$ edges. 
Observe that $\S$ is a partition of $V(X)$, so we can also consider the quotient graph $X_\S$ which is precisely
the underlying graph of $X_\rho$. \bigskip

\noindent
{\bf Remark.} Note that if $G$ is as in Lemma~\ref{lem:blocksemi} and $\varphi \in G$ is $(m,q)$-semiregular,
then the subgroup $\la \varphi \ra$ has $m$ orbits of length $q$, and so Lemma~\ref{lem:blocksemi} applies.\bigskip

For the sake of completeness we state the following classical result which will be
used throughout the paper.

\begin{proposition}{\rm \cite[Theorem~3.4]{W64}}
\label{pro:wielandt}
	Let $p$ be a prime and let $P$ be a Sylow $p$-subgroup of a permutation group $G$ acting on a set $\Omega$.
	Let $\omega \in \Omega$. If $p^m$ divides the length of the $G$-orbit containing $\omega$,
	then $p^m$ also divides the length of the $P$-orbit containing $\omega$. 
\end{proposition}

The following proposition is a generalization of \cite[Theorem~3.4]{DM81}.
 
 
 \begin{proposition}
 \label{pro:semireg}
 Let $X$ be a vertex-transitive graph of order $mp$, where $m < p$, $p$ a prime, and let $G \leq \Aut X$ be
 a transitive subgroup of automorphisms of $X$. Then there exists some $(m,p)$-semiregular automorphism
 $\rho$ of $X$, such that $\rho \in G$.
 \end{proposition}
 
 \begin{proof}
 Since $G$ is transitive on $V(X)$ and $X$ is of order $mp$, the order $|G|$ of $G$ is divisible by $p$.
 Let $P$ be a Sylow $p$-subgroup of $G$. Since the length $l$ of an orbit of $P$ divides its order $|P|$, 
 it can either be $1$ or $p$ (recall that $m < p$). By Proposition~\ref{pro:wielandt}, $p$ divides $l$ and thus
 $l = p$. Therefore $P$ has exactly $m$ orbits of length $p$. Following the proof of \cite[Theorem~3.4]{DM81}
 one can now show that there exists some $\rho \in P$ such that $\rho$ is $(m,p)$-semiregular.
 \end{proof}
 
 The following lemma can be deduced from \cite[Lemma~2]{D95}.
 \begin{lemma}
 \label{le:semiregnormal}
 Let $X$ be a vertex-transitive graph of order $mq$, where $q$ is a prime, 
 let $G$ be an imprimitive subgroup of
 automorphisms of $X$ and let $N$ be a normal subgroup of $G$ with orbits of length $q$. Then $X$ has an 
 $(m,q)$-semiregular automorphism whose orbits coincide with the orbits of $N$.
 \end{lemma}

We now introduce the following notion of a lift of a path in a graph with a semiregular automorphism.
Let $X$ be a graph that admits an $(m,n)$-semiregular automorphism $\rho$. Let $\S = \{S_1, S_2, \ldots , S_m\}$
be the set of orbits of $\rho$, let $X_\S$ be the corresponding quotient graph and let $\p : X \to X_\S$ 
be the corresponding projection.
Let $W = S_{i_1}S_{i_2}\ldots S_{i_k}$ be a path in $X_\S$. We let the {\em lift of the path} $W$ be
the set of all paths of $X$ whose images under $\p$ are $W$. The following lemma is straightforward and is
just a reformulation of \cite[Lemma~5]{MP82}.

\begin{lemma}
\label{lem:cyclelift}
Let $X$ be a graph admitting an $(m,p)$-semiregular automorphism $\rho$, where $p$ is a prime. 
Let $C$ be a cycle of length $k$ in the quotient graph $X_\S$, where $\S$ is the set of orbits of $\rho$. 
Then, the lift of $C$ either contains a cycle of length $kp$ or it consists of $p$ disjoint $k$-cycles.
In the latter case we have $d(S,S') = 1$ for every edge $SS'$ of $C$.
\end{lemma}

A path of $X$ which meets each of the vertices of $X$ is called a {\em Hamilton path} of $X$.
A Hamilton cycle is defined in a similar way. The following classical result, due to Jackson \cite{J78},
giving a sufficient condition for the existence
of Hamilton cycles in $2$-connected regular graphs will be used throughout this paper (Note that 
every connected vertex-transitive graph is $2$-connected).

\begin{proposition}
\label{pro:jack}
{\rm \cite[Theorem~6]{J78}}
Every $2$-connected regular graph of order $n$ and valency at least $n/3$
contains a Hamilton cycle.
\end{proposition}

The next result may be extracted from \cite[Theorem~2.10]{DMMN}.

\begin{theorem}
\label{the:simple}
Let $G$ be a transitive permutation group of degree $6p$, $p \ge 5$ a prime, 
with an imprimitivity block system $\B$ formed
by a (proper, intransitive) minimal normal subgroup $N$ of $G$. Then  $N^B$ is simple
for all blocks $B\in\B$. 
\end{theorem}

We let $\ZZ_n =\{0, 1,\ldots, n-1\}$ denote the ring of integers modulo
$n$, and we let $\ZZ_n^*$ be the multiplicative group of the units of $\ZZ_n$.


In the subsequent sections some of the graphs will be represented in the Frucht's notation~\cite{RF70}.
For the sake of completeness we include the definition.
Let $X$ be a connected vertex-transitive graph of order $mn$ admitting an $(m,n)$-semiregular 
automorphism $\rho$. Let $\S = \{S_i\ |\ i \in \ZZ_m\}$ be the set of orbits of $\rho$.
Denote the vertices of $X$ by $v_i^j$, where $i \in \ZZ_m$ and $j \in \ZZ_n$, in such a way that
$S_i = \{v_i^j\ |\ j \in \ZZ_n\}$ with $v_i^j = v_i^0\rho^j$. 
Then $X$ may be represented by the notation of Frucht~\cite{RF70}
emphasizing the $m$ orbits of $\rho$ in the following way. 
The  $m$ orbits of $\rho$ are  represented by $m$ circles.
The symbol $n/R$, where $R \subseteq\ZZ_n\setminus \{0\}$, inside a circle corresponding to   
the orbit $S_i$ indicates that for each $j\in \ZZ_n$,
the vertex $v_i^j$ is adjacent to all the vertices $v_i^{j+r}$, where $r \in R$. 
When $X\la S_i \ra$ is an independent set of vertices we simply write $n$ inside its circle.
Finally, an arrow pointing from  the circle representing the orbit $S_i$ to the circle representing the
orbit $S_k$, $k \ne i$, labeled by the set $T \subseteq \ZZ_n$ indicates that
for each $j\in\ZZ_n$, the vertex $v_i^j\in S_i$ is adjacent to all the vertices $v_{k}^{j+t}$,
where $t \in T$.
An example illustrating this notation is given in Figure~\ref{fig:petersen}.


\section{Genuinely imprimitive graphs}
\label{sec:three}
\indent

Throughout this section let $X$ be a connected genuinely imprimitive graph of order $6p$, $p > 3$ a prime,
admitting an imprimitive subgroup $G$ of $\Aut X$ with a  nontransitive 
minimal normal subgroup $N \triangleleft G$.
Let the set of orbits of $N$ (and thus blocks for $G$) be denoted by $\B$.

The task of showing that $X$ has a Hamilton path is divided into 
six different cases depending on the size of the blocks in $\B$.
Each of them is covered by a separate lemma 
(see Lemmas~\ref{lem:2}, \ref{lem:3}, \ref{lem:p}, \ref{lem:6}, \ref{lem:2p} and \ref{lem:3p}).
If the size of blocks equals to $p$ or $6$ we in fact show that $X$ contains a Hamilton cycle.

\begin{lemma}
\label{lem:2}
If the size of blocks in $\B$ is $2$ then $X$ 
has a Hamilton path. 
\end{lemma}

\begin{proof}
Since $X_\B$ is a connected vertex-transitive graph of order $3p$ it has a Hamilton cycle $C$. 
By Lemma~\ref{le:semiregnormal}, $X$ has a $(3p,2)$-semiregular automorphism whose  set of orbits
equals $\B$.
Thus, by Lemma~\ref{lem:cyclelift}, the lift of $C$ either contains a Hamilton cycle of $X$ or it contains a
disjoint union of two cycles of length $3p$. Since $X$ is connected a Hamilton path exists in $X$.
\end{proof}

The following auxiliary lemma will be used in the proof of Lemma~\ref{lem:3}.

\begin{lemma}
\label{lem:3Pet}
If the size of blocks in $\B$ is $3$ and the quotient graph $X_\B$ is isomorphic to the Petersen
graph then $X$ has a Hamilton path.
\end{lemma}

\begin{proof}
Note that in this case $p = 5$.
By Lemma~\ref{le:semiregnormal} there exists a $(10,3)$-semiregular automorphism $\varphi$ of $X$
whose orbit set equals $\B$. Suppose there exist two disjoint $5$-cycles in $X_\B$ whose lifts both
contain a $15$-cycle. Then the connectedness of $X$ implies that $X$ has a Hamilton path. We can thus
assume that no two such $5$-cycles exist in $X_\B$. We claim that this implies that for any two
adjacent orbits $B, B' \in \B$ of $\varphi$ we have $d(B,B') = 1$. Suppose this is not the case. It is easy to see that
we then have two disjoint $5$-cycles in $X_\B$ such that each of them contains an edge corresponding to a
multiedge in $X_\varphi$. But then Lemma~\ref{lem:cyclelift} implies that the lifts of both of these two
$5$-cycles contain $15$-cycles, a contradiction.

Note that in the case when $X\la B \ra$ is not an independent set for some (and thus all) 
$B \in \B$ a Hamilton path exists in $X$. We can thus assume
that $X\la B \ra$ is an independent set for all $B \in \B$.
Let $\bar{G}$ denote the permutation group corresponding to the natural action of $G$ on $X_\B$. Since
the only transitive subgroups of the automorphism group of the Petersen graph are $S_5$, $A_5$ and 
$AGL(1,5)$, the fact that $\bar{G}$ is transitive implies, that a subgroup $H$ of $\bar G$, which is
isomorphic to $\AGL(1,5)$ or to $A_5$, exists. As we demonstrate below, each of these two cases lead to a contradiction, 
which shows that $X$ has a Hamilton path, as required.

Suppose first that $H \cong \AGL(1,5)$.
Then there exist two disjoint $5$-cycles of $X_\B$ interchanged by some element of $H$. 
The lift of each of them
is thus a union of $3$ disjoint $5$-cycles. Hence, we can assume that the Frucht's notation of $X$ is as
in Figure~\ref{fig:petersen}. In view of our assumptions we have
$$ a = c\ \mathrm{or}\ d=e, \quad b = d\ \mathrm{or}\ a = e, \quad c = e\ \mathrm{or}\ a = b, \quad
	a = d\ \mathrm{or}\ b = c, \quad b = e\ \mathrm{or}\ c=d.$$
As $X$ is connected, we cannot have $a=b=c=d=e$. With no loss of generality assume that $a \neq b$, and so  
$c = e$. 
Suppose first that $a = d$. Then $b \neq d$, and so $d = a = e = c$. The reader may check that then
the vertices of $B_1$ are contained on precisely two $5$-cycles, whereas the vertices of $B_0$ are contained on
precisely four $5$-cycles which is impossible in view of vertex-transitivity of $X$. Suppose then that 
$a \neq d$. Therefore, $b = c$ and thus also $d = e = c = b$. As above a contradiction to vertex-transitivity of $X$ is 
obtained.

\begin{figure}[htb]
\begin{minipage}[b]{0.5\linewidth}
\centering
\includegraphics[width=0.8\linewidth]{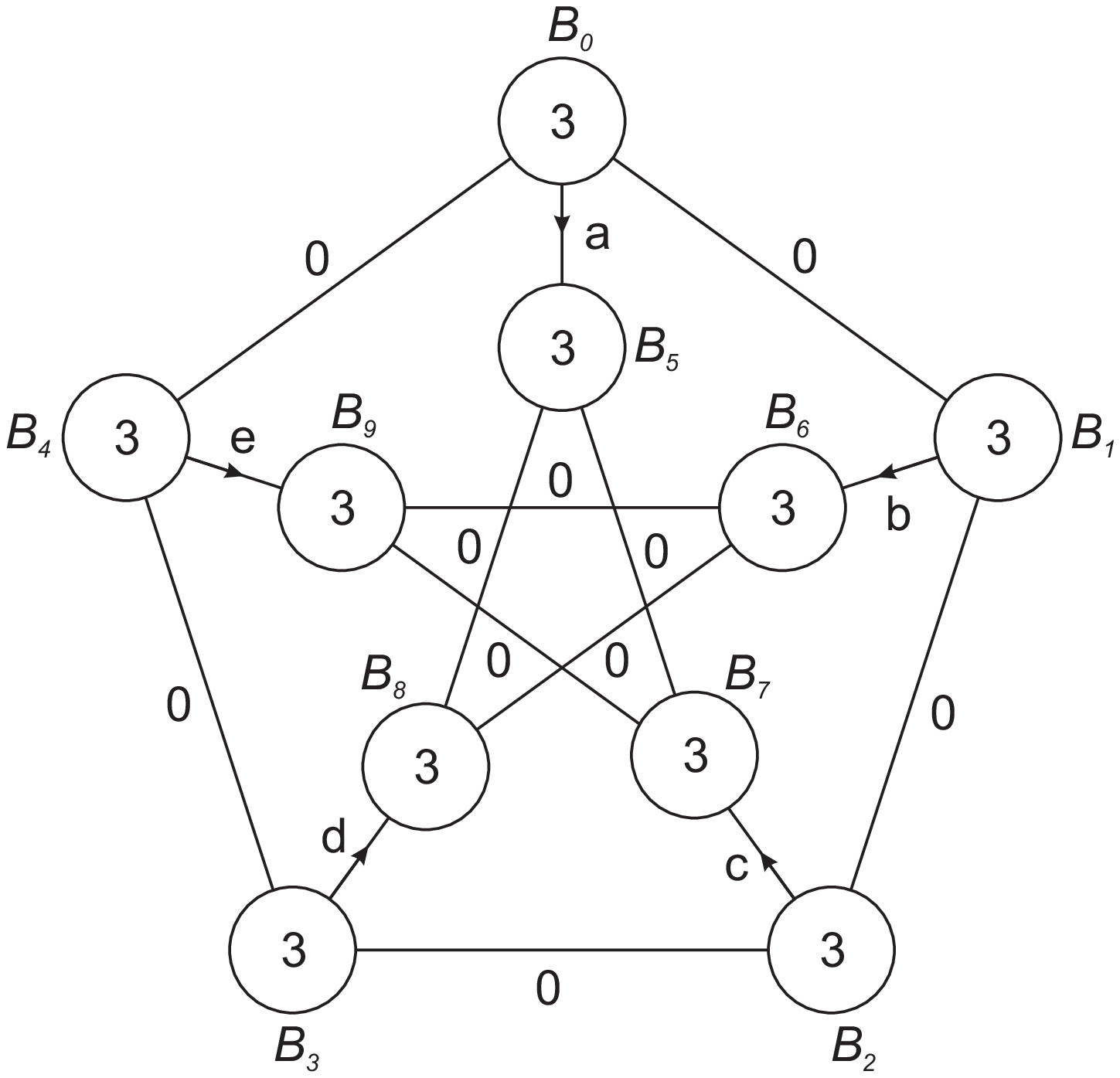}
\caption{\label{fig:petersen}  The case $H = \AGL(1,5)$.}
\end{minipage}\hspace{0.5cm}
\begin{minipage}[b]{0.5\linewidth}
\centering
\includegraphics[width=0.8\linewidth]{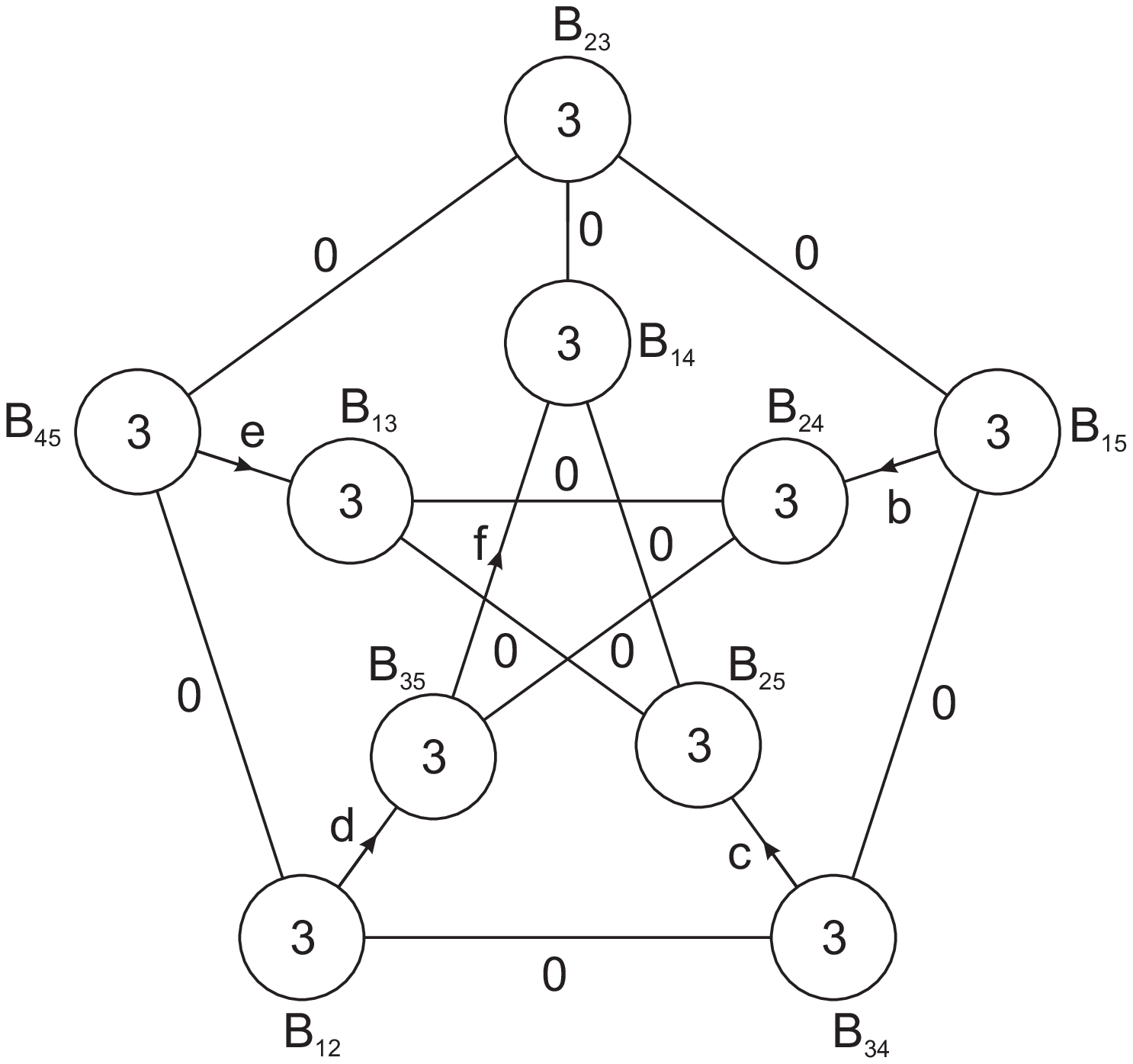}
\caption{\label{fig:petersenA5}  The case $H = A_5$. }
\end{minipage}
\end{figure}

Suppose now that $H \cong A_5$. We can assume that the Frucht's notation of $X$ is as in
Figure~\ref{fig:petersenA5}, where the group $H$ acts on $X_\B$ in the obvious way.
In view of the action of an automorphism of $H$ whose action on $X_\B$ corresponds to $(23)(45)$, we have
$e = 0$. Furthermore, the element of $H$ corresponding to $(12)(45)$ forces $d = 0$. 
Continuing in this way we find that
$c-f=0$, $b-c=0$ and $b+f=0$, which forces $b = c = f = 0$. However, this contradicts the connectedness of $X$
and the proof is completed.
\end{proof}

An $n$-{\em bicirculant} is a graph with a $(2,n)$-semiregular automorphism.
Every $n$-bicirculant $X$ can be represented 
by a triple of subsets of $\ZZ_n$
in the following way. Let $\varphi$ be a $(2,n)$-semiregular automorphism of $X$,
let $U$ and $W$ be the two orbits of $\varphi$,
and let $u \in U$ and $w \in W$. 
Let $S = \{ s \in \ZZ_n\ |\ u \sim u\varphi^s \}$
be the symbol of the $n$-circulant induced on $U$ and 
let $R$ be the symbol of the $n$-circulant induced on $W$ (relative to $\varphi$).
Moreover, let $T =\{ t \in \ZZ_n\ |\ u \sim w\varphi^t \}$.
The ordered triple $[S,R,T]$ is the {\em symbol of} $X$ {\em relative to} 
$(\varphi,u,w)$. Note that $S=-S$ and $R=-R$ are symmetric, that is, inverse-closed
subsets of $\ZZ_n$, and are independent 
of the particular choice of vertices $u$ and $w$.

In the rest of this section the well known wreath and Cartesian products of graphs will be encountered.
To fix the notation, we include the definitions.
For two graphs $X$ and $Y$ let $X\wr Y$ denote the {\em wreath product} of $X$ by $Y$,
that is, the graph with 
vertex set $V(X)\times V(Y)$ with two vertices $(a,u)$ and $(b,v)$
adjacent in  $X\wr Y$  if and only if either  $ab\in E(X)$ or  $a=b$ and $uv\in E(Y)$.
Note that the wreath product is sometimes refered to as the {\em lexicographic product}.
The {\em Cartesian product} $X\square Y$ of graphs $X$ and $Y$ 
is the graph with vertex set $V(X)\times V(Y)$, where two vertices $(a,u)$ and $(b,v)$ are adjacent
in  $X\square Y$  if and only if either  $ab\in E(X)$ and  $u=v$, or $a=b$ and $uv\in E(Y)$. 

\begin{lemma}
\label{lem:3}
If the size of blocks in $\B$ is $3$ then $X$ 
has a Hamilton path. 
\end{lemma}

\begin{proof}
By Lemma~\ref{le:semiregnormal} there exists a $(2p,3)$-semiregular automorphism $\varphi$ of $X$ whose
orbit set coincides with $\B$. If the quotient graph $X_\B$ is isomorphic to the Petersen graph, then
Lemma~\ref{lem:3Pet} applies.
We can thus assume that $X_\B$ is not isomorphic to the Petersen graph. Therefore, $X_\B$ has a
Hamilton cycle $C = B_0B_1 \ldots B_{2p-1}B_0$.
In view of Lemma~\ref{lem:cyclelift} we can assume that the lift of $C$ consists of 
three disjoint $2p$-cycles. So $d(B_i,B_{i+1}) = 1$ for all $i \in \ZZ_{2p}$.
Therefore, we can label the vertices of $X$ by $\{u_i^j\ |\  i \in \ZZ_{2p},\ j \in \ZZ_3\}$ in such a way
that $B_i = \{u_i^j\ |\ j \in \ZZ_3\}$ and that $u_i^ju_{i+1}^j$ is an edge of $X$ for every
$i \in \ZZ_{2p}$ and $j \in \ZZ_3$. Moreover, we can assume that
$X\la B\ra = 3K_1$ for all $B \in \B$  (otherwise $X$ contains a subgraph isomorphic to the 
Cartesian product $C_{2p}\square K_3$ which clearly has a Hamilton cycle).

There exists some $\psi \in N$ such that $\psi|_{B_0} = (u_0^0u_0^1u_0^2)$.
By the above assumptions it is clear that $\psi|_{B_i} = (u_i^0u_i^1u_i^2)$ for all $i \in \ZZ_{2p}$.
Therefore, we can assume that the automorphism $\varphi$ is in $N$.
Note also that $N$ acts faithfully on each of its orbits $B \in \B$ and thus either $N \cong \ZZ_3$ or $N \cong S_3$.
However, the latter case cannot occur, for then the Sylow $3$-subgroup of $N$ is normal in $G$, 
contradicting the minimality of $N$.

By Proposition~\ref{pro:semireg} a $(6,p)$-semiregular automorphism of $X$ exists if $p > 5$.
We now show that such an automorphism exists also if $p = 5$.
Suppose then that $X$ is of order $30$. 
Let $P \leq G$ be a Sylow $5$-subgroup of $G$. By Proposition~\ref{pro:wielandt} the lengths of its orbits
are divisible by $5$. Therefore, $P$ either has $6$ orbits of length $5$ or one 
orbit of length $25$ and one orbit of length $5$. However, a similar argument as in the proof
of Lemma~\ref{lem:blocksemi} shows that the latter case is impossible. 
So $P$ has $6$ orbits of length $5$. By Lemma~\ref{lem:blocksemi} it follows
that the group $P$ has two orbits of length $5$ in its natural action on $X_\B$. Thus an element $\psi \in P$ 
of order $5$ is either $(6,5)$-semiregular or it has $3$ orbits of length $5$ and $15$ fixed points. 
In the latter case there exists some other element $\vartheta \in P$ such that none of the above 
$15$ fixed points of $\psi$
is fixed by $\vartheta$. Hence either $\vartheta$ or $\vartheta\psi$ is $(6,5)$-semiregular. This proves
our claim that a $(6,p)$-semiregular automorphism of $X$ always exists. Let us denote it by $\rho$.

We claim that $\rho$ and $\varphi$ commute. 
Namely, since $N \cong \ZZ_3$, we have that $\rho^{-1}\varphi\rho$ is equal either to $\varphi$ or $\varphi^{-1}$.
But $p$ is odd, so $\rho^{-1}\varphi\rho = \varphi^{-1}$ would imply $\rho^{-p}\varphi\rho^p = \varphi^{-1}$,
which is clearly impossible as $\rho^p = 1$. 
Thus, $\varphi\rho = \rho\varphi$. Moreover, this element is of order $3p$ 
and has precisely two orbits of length $3p$ which implies that $X$ is a bicirculant. 
Let $[S,R,T]$ be one of its symbols corresponding to $\varphi\rho$, such that $0 \in T$. 
If there exists some $a \in T$ for which $\la a \ra = \ZZ_{3p}$, where $\la a \ra$ is the additive subgroup
of $\ZZ_{3p}$ generated by $a$, then $X$ has a Hamilton cycle.  
Moreover, if $T$ contains an element
of order $p$ and an element of order $3$, then their difference generates $\ZZ_{3p}$, and so $X$ has a Hamilton
cycle. We can therefore assume that $\la T \setminus \{0\} \ra$ is either empty or it is one of 
$\la 3 \ra$ or $\la p\ra$. 

As $X\la B \ra$ is an independent set for each $B \in \B$, there is no element of order $3$ in $S$ or in $R$.
If $\la S \ra = \ZZ_{3p}$ and $\la R \ra = \ZZ_{3p}$, then the subgraphs induced on each of the orbits of 
$\varphi\rho$ are connected vertex-transitive graphs of order $3p$,
and so they both contain a Hamilton cycle. Clearly, $X$  has a Hamilton path in this case.
With no loss of generality we can thus assume that $\la S \ra \neq \ZZ_{3p}$. This implies
that $S = \emptyset$ or $\la S \ra = \la 3 \ra$.
Suppose first that $S = \emptyset$. Then regularity of $X$ implies 
$R = \emptyset$ as well. By the above remarks on $T$, $X$ is not connected, a contradiction.
Therefore, $\la S \ra = \la 3 \ra$. As $X$ is regular, we have that $|S| = |R|$, and so
either $\la R \ra = \la 3 \ra$ or $\la R \ra = \ZZ_{3p}$. In the former case 
the subgraph induced on each of the orbits of $\rho$ contains a $p$-cycle. Moreover, the facts that
$\la T \ra \neq \ZZ_{3p}$ and $X$ is connected imply, that there exists some $a \in T$ of order $3$,
and so $a$ and $0$ give rise to a $6$-cycle of $X_\rho$. Therefore, $X$ has a Hamilton path in this case.
We are left with the possibility $\la R \ra = \ZZ_{3p}$. In view of
the fact that no element of order $3$ exists in $R$, some $a \in R$ such that $\la a \ra = \ZZ_{3p}$ exists.
We can assume that $a = 1$ (otherwise take $(\varphi\rho)^a$ instead of $\varphi\rho$). Since $\la S \ra = \la 3 \ra$, we
have $3k \in S$ for some $k \in \{1,2,\ldots , p-1\}$. Thus $X$ contains a subgraph isomorphic to the
generalized Petersen graph $GP(3p, 3k)$ which
has a Hamilton cycle (see \cite{BA83}).
\end{proof}

\begin{lemma}
\label{lem:p}
If the size of blocks in $\B$ is $p$ then $X$ 
has a Hamilton cycle.
\end{lemma}

\begin{proof}
The quotient graph $X_\B$ is a connected vertex-transitive graph on $6$ vertices. 
By Lemma~\ref{lem:blocksemi} the blocks of $\B$ coincide with the orbits of some $(6,p)$-semiregular 
automorphism $\rho \in G$ of $X$,
which exists by Lemma~\ref{le:semiregnormal}. Let $\S=\{S_i\mid i\in \ZZ_6\}$ denote the set of orbits
of $\rho$ and denote the vertices of each $S_i$ with $u_i^j$, $j\in\ZZ_p$, where $u_i^j\rho = u_i^{j+1}$. 
The quotient graph
$X_\S=X_\B$ is isomorphic to one of the following five graphs: $C_6$, $K_3\square K_2$, $K_{3,3}$,
$K_3\wr 2K_1$ or $K_6$ (these are the only connected vertex-transitive graphs on six vertices).
It is easy to see that in all these cases for any edge $e$ of $X_\S$ there exists a Hamilton cycle 
of $X_\S$ containing $e$. Hence, by Lemma~\ref{lem:cyclelift}, 
we may assume that no multiedge exists in $X_\rho$.  
Moreover, we may label the orbits of $\rho$ in such a way that $S_i\sim S_{i+1}$ for every $i\in\ZZ_6$.
If there exists a Hamilton cycle of $X_\S$ whose lift contains a Hamilton cycle of $X$, there is
nothing to prove. Therefore, we can assume that no such Hamilton cycle of $X_\S$ exists.
Consequently, we may assume that $u_i^j\sim u_{i+1}^j$, $i\in\ZZ_6$ and $j\in\ZZ_p$.
Note also that we can assume that $X\la S_i \ra = pK_1$ for all $i \in \ZZ_6$. 
Namely, if the subgraphs $X\la S_i \ra$ are of valency $2$, then a Hamilton cycle of $X$ exists by 
\cite[Theorem~3.9]{BA89}, and if the subgraphs $X\la S_i \ra$ are of valency at least $4$, 
then \cite[Theorem~4]{CQ81} implies that each of $X\la S_i \ra$ is Hamilton-connected
(that is, there exists a Hamilton path of $X\la S_i \ra$ connecting any two vertices), 
and so a Hamilton cycle of $X$ clearly exists. 

We distinguish five different cases depending on which of the five connected vertex-transitive graphs of order $6$
the quotient graph $X_\S$ is isomorphic to.

If $X_\S\cong C_6$ then $S_iS_{i+1}$, where $i\in\ZZ_6$, are the only edges of $X_\S$, and so
$X$ is not connected, a contradiction.

Suppose that $X_\S\cong K_3\square 2K_1$. Then we may assume that in addition to
the edges $S_iS_{i+1}$, also $S_0S_4, S_1S_3, S_2S_5 \in E(X_\S)$.
Therefore, 
$$
E(X)=\{u_i^ju_{i+1}^j\mid i\in\ZZ_6, j\in\ZZ_p\}\cup \{u_0^ju_4^{j+r_0}, u_1^ju_3^{j+r_1}, u_2^ju_5^{j+r_2} \mid j\in \ZZ_p\},
$$
where $r_0,r_1,r_2\in\ZZ_p$. Since $S_0S_4S_3S_1S_2S_5S_0$ 
and $S_0S_1S_3S_2S_5S_4S_0$ are  Hamilton cycles of $X_\S$, 
Lemma~\ref{lem:cyclelift} implies that $r_0-r_1+r_2=0$ and $r_0-r_2-r_1=0$. 
Subtracting one of the equations from the other we get that $r_2=0$, and so $r_0=r_1$. In view of the connectedness
of $X$, we have $r_0=r_1\ne 0$. Then
$$
u_0^0u_4^{r_0}u_5^{r_0}u_0^{r_0}u_4^{2r_0}\cdots u_0^{-r_0}u_4^0u_5^0u_2^0u_3^0u_1^{-r_0}u_2^{-r_0}u_3^{-r_0}u_1^{-2r_0}\cdots u_2^{r_0}u_3^{r_0}u_1^0u_0^0
$$
is a Hamilton cycle of $X$.

Suppose next that $X_\S\cong K_{3,3}$. 
Hence we may assume that adjacencies in $X_\S$ are 
$S_i\sim S_{i+1}$ and $S_i\sim S_{i+3}$, where $i\in\ZZ_6$. 
This implies that 
$E(X)=\{ u_i^ju_{i+1}^j, u_i^ju_{i+3}^{j+r_i}\mid i\in \ZZ_6,  j, r_i\in\ZZ_p\}$, where  $r_i=-r_{i+3}$. 
Since  $S_0S_3S_2S_1S_4S_5S_0$, $S_0S_3S_4S_5S_2S_1S_0$ and
$S_0S_3S_2S_5S_4S_1S_0$ are Hamilton cycles of $X_\S$, Lemma~\ref{lem:cyclelift} implies
that $r_0+r_1=0$, $r_0+r_5=r_0-r_2=0$ and  $r_0+r_2+r_4=r_0+r_2-r_1=0$. 
As $p \geq 5$, combining these equations we get that $r_i=0$ for every $i\in\ZZ_p$, which contradicts
the fact that $X$ is connected.   

The remaining two cases ($X_\S = K_3\wr 2K_1$ and $X_\S = K_6$) are dealt with in a similar manner. We leave the
details to the reader.
\end{proof}

\noindent
{\bf Remark.} In the above proof a Hamilton cycle was shown to exist in $X$ using the following idea.
When considering the possible arrangements of the edges of $X$, where the quotient graph $X_\S$ has been given,
the key factors are the connectedness of $X$ and Lemma~\ref{lem:cyclelift}. 
This way we find that either a Hamilton cycle of $X_\S$ whose lift
contains a Hamilton cycle of $X$ exists, or the structure of the edges of $X$ is completely determined in which
case a Hamilton cycle of $X$ is easily found. The same approach will be used throughout this paper. The
technical details will be left to the reader.\bigskip\bigskip

\setcounter{case}{0}
\begin{lemma}
\label{lem:6}
If the size of blocks in $\B$ is $6$ then $X$ 
has a Hamilton cycle. 
\end{lemma}

\begin{proof}
Note that $X_\B$ is a connected $p$-circulant so it has a Hamilton cycle.
Theorem~\ref{the:simple} implies, that  $N^B$ is simple of degree $6$ for every $B\in\B$. 
The only two transitive simple groups of degree $6$ 
up to permutation isomorphism are
the alternating group $A_6$ and its subgroup isomorphic to $A_5$ (see \cite{DM96}). They are both
doubly transitive. Thus 
the   subgraphs  $X\la B\ra$, $B \in \B$, are either all isomorphic to $K_6$  or they are
all isomorphic to $6K_1$. 

Suppose first that $X\la B\ra$ is isomorphic to $K_6$ for all $B \in \B$.  Then $X\la B\ra$  
is Hamilton connected for every $B \in \B$, and so a Hamilton cycle of $X$ clearly exists.

Suppose now that $X\la B\ra =6K_1$ for all $B \in \B$. 
Every simple subgroup of $A_6$ of order $60$ is permutation isomorphic
to $H = \la (1\,2\,3\,4\,5), (1\, 2)(4\,6)\ra$ (see for example \cite[Table~2.1]{DM96}). 
Thus for any $B \in \B$ and any vertex $v \in B$ we
have some $\alpha \in N^B$ fixing $v$ and cyclically permuting the other five vertices of $B$.
We claim that for any two adjacent blocks $B, B' \in \B$ the
graph $[B,B']$ is isomorphic to $K_{6,6}$, to $K_{6,6} - 6K_2$ or to $6K_2$. Namely, suppose that
a vertex $u \in B$ has at least two neighbors, say $v_1$ and $v_2$, in $B'$. By the above remarks there exists
an automorphism $\alpha \in N$ fixing $u$ and permuting the other five vertices of $B$. We distinguish 
two different cases depending on the order $d$ of $\alpha|_{B'}$. 

\begin{case}\end{case} $d = 5$. Then $\alpha|_{B'}$ also fixes a vertex $v$ of $B'$ and cyclically permutes the
other five vertices of $B'$. With no loss of generality assume $v \neq v_1$.
Applying $\alpha$ to the edge $uv_1$ we get that the valency of
$u$ in $[B,B']$ is either $5$ or $6$, depending on whether
$u$ is adjacent to $v$ or not. Since $\B$ is the set of orbits of $N$, a simple counting argument shows that
the subgraph $[B,B']$ is isomorphic either to $K_{6,6}$ or to $K_{6,6} - 6K_2$ as claimed. 

\begin{case}\end{case} $d \neq 5$. With no loss of generality we can assume that $d = 1$ (otherwise take an 
appropriate power of $\alpha$). Since $u$ has a neighbor in $B'$, every vertex of $B$ has a neighbor in $B'$.
Let $u' \in B$, $u' \neq u$, have a neighbor $v$ in $B'$. Applying $\alpha$ to the edge $u'v$ we get that
$v$ is adjacent to all the vertices of $B$ except possibly $u$. Thus
 $[B,B']$ is isomorphic either to $K_{6,6}$ or to $K_{6,6} - 6K_2$, which completes the proof of our claim.\smallskip

Now let $B \in \B$. We claim that  there exists a block $B'$, adjacent to $B$, such that $[B,B']$
is not isomorphic to $6K_2$. Namely, if this is not the case, then a contradiction to the connectedness of $X$ is
obtained by an argument similar to the one of the above two paragraphs.
Since $G$ acts transitively on $X$, there
exists an element $\psi \in G$ cyclically permuting the $p$ blocks of $\B$. 
With no loss of generality we can assume that $B' = B\psi$ (otherwise take
an appropriate power of $\psi$). It follows that $B\psi^i \sim B\psi^{i+1}$ for all $i \in \ZZ_p$.
It is now evident that $X$ has a Hamilton cycle. 
\end{proof}


\begin{lemma}
\label{lem:2p}
If the size of blocks in $\B$ is $2p$ then $X$ 
has a Hamilton path. 
\end{lemma}

\begin{proof}
Note that $X_\B=K_3$ and that the group $G$ acts edge transitively in its natural action on $X_\B$.
Let $\B = \{B_i\ |\ i \in \ZZ_3\}$. Let $P \leq G$ be some Sylow $p$-subgroup of $G$.
In view of Proposition~\ref{pro:wielandt} and the fact that $G$ has $3$ blocks of size $2p$,
$P$ has $6$ orbits of length $p$. Denote them by $\S = \{S_i\ |\ i \in \ZZ_6\}$. 
By Lemma~\ref{lem:blocksemi} each block in $\B$ is a union of two orbits of $P$.
With no loss of generality we can assume that $B_0 = S_0 \cup S_1$, $B_1 = S_2 \cup S_3$ and $B_2 = S_4 \cup S_5$.
 
By Proposition~\ref{pro:semireg}, there exists a $(6,p)$-semiregular automorphism
$\rho$ of $X$ such that $\rho \in G$ whenever $p > 5$. We show that we can assume such an element to exist
even if $p = 5$. To this end suppose that $p = 5$ and that $X$ does not contain a Hamilton path. 
In view of Proposition~\ref{pro:jack}
the valency of $X$ is at most $9$. Let $\rho \in G$ be an element of order $5$, whose action
on $B_0$ is $(2,5)$-semiregular (which exists by Proposition~\ref{pro:semireg}). 
With no loss of generality assume that $\rho \in P$. The two orbits of $\rho$ in $B_0$ thus coincide with $S_0$
and $S_1$. If $\rho$ is not $(6,5)$-semiregular, then we can assume that it fixes some vertex $u \in S_2$.
Since $X_\B = K_3$, the vertex $u$ has a neighbor in $B_0$ and thus its valency in $[B_0,B_1]$ is at least $5$.
As $[B_0,B_1]$ is regular and $G$ acts edge-transitively on $X_\B$, the valency of $u$ in 
$[B_1,B_2]$ is at least $5$ as well, contradicting the fact that $u$ has valency at most $9$.
Thus $\rho$ is $(6,5)$-semiregular, as required.

We can clearly assume that the orbit set of $\rho$ is $\S$. 
In view of regularity of the bipartite graphs $[B,B']$, $B,B' \in \B$, the subgraph $\bar{X}_\S$
of $X_\S$, which is obtained from $X_\S$ by deleting the edges $S_0S_1$, $S_2S_3$, $S_4S_5$ 
(if they exist), is
clearly one of the graphs $Y_i$, $i \in \{0,1,2,3,4\}$ of Figure~\ref{fig:rho2p}.
However, for each of the graphs $Y_i$, $i \geq 1$, the following holds: if there exists a multiedge of $X_\rho$, 
then there exists a Hamilton cycle of $X_\S$ which contains an edge
corresponding to a multiedge of $X_\rho$. By Lemma~\ref{lem:cyclelift} we can thus assume that
no multiedge exists in $X_\rho$, except possibly if $\bar{X}_\S = Y_0$. In view of the regularity of $X$
the graphs $Y_3$ and $Y_4$ are then not possible.

\begin{figure}[htbp]
\begin{center}
\includegraphics[width=0.90\hsize]{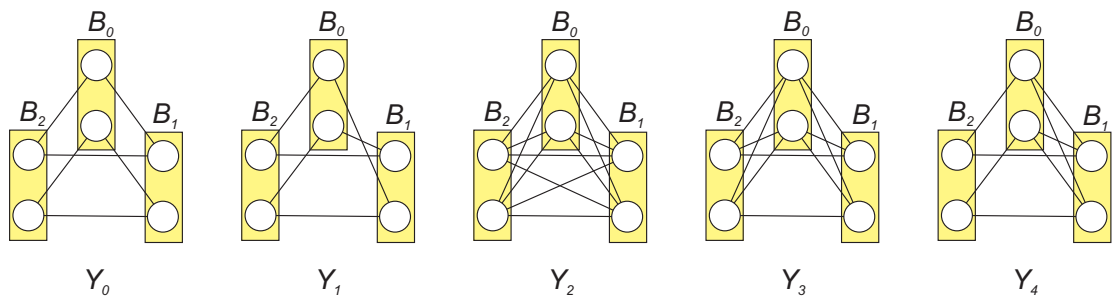}
\caption{\label{fig:rho2p} All possibilities for the subgraph  $\bar{X}_\S$ of $X_\S$
when $G$ has a complete block system $\B$ of three blocks of size $2p$.}
\end{center}
\end{figure}

If $X\la B_0 \ra$ is a connected graph, then for each of its vertices there exists a Hamilton path of $X\la B_0 \ra$
starting at that vertex, so $X$ clearly has 
a Hamilton path in this case. We can thus assume that $X\la B_0 \ra$ is not connected. As it is a vertex-transitive
graph, it is isomorphic to $2pK_1$, to $pK_2$ or it is a disjoint union of two isomorphic connected $p$-circulants.
We consider each of the three cases separately.

\setcounter{case}{0}
\begin{case} $X\la B_0 \ra \cong 2pK_1$. \end{case}
As $X$ is connected, the quotient graph $X_\S = \bar{X}_\S$ is one of
$Y_1$ or $Y_2$. If $X \cong Y_1$, then connectedness of $X$ and Lemma~\ref{lem:cyclelift} imply that
the lift of $Y_1$ contains a Hamilton cycle of $X$. 
It is easy to see that if $X_\S \cong Y_2$, the connectedness of $X$ forces 
some Hamilton cycle of $X_\B$, whose lift contains a Hamilton cycle of $X$, to exist. 
We leave the details to the reader.

\begin{case} $X\la B_0 \ra \cong pK_2$. \end{case} It is clear that $[S_0,S_1] \cong pK_2$. 
Suppose first that $\bar{X}_\S \cong Y_0$. In this case every edge of $X_\S$ is contained on some
Hamilton cycle of $X_\S$, and so Lemma~\ref{lem:cyclelift} 
implies that we can assume that no multiedge exists in $X_\rho$. 
If there exists a Hamilton cycle of $X_\S$ whose lift contains a Hamilton cycle of $X$, we are done. If not,
the connectedness of $X$ implies that $X \cong C_{3p} \square K_2$, and so $X$ contains a Hamilton cycle.
In the case when $\bar{X}_\S$ is isomorphic to one of $Y_1$ and $Y_2$ one can easily see that
the connectedness of $X$ forces some Hamilton cycle of $X_\S$, whose lift contains a Hamilton cycle of $X$, to exist.
The details are left to the reader. 

\begin{case} $X\la B_0 \ra$ \end{case} is isomorphic to a disjoint union of two isomorphic connected
$p$-circulants. In view of connectedness of $X$ the quotient graph $X_\S = \bar{X}_\S$ is one of
$Y_1$ or $Y_2$ and so it has a Hamilton cycle. 
As the six $p$-circulants are precisely the graphs $X\la S_i \ra$, where $i \in \ZZ_6$, 
a Hamilton path exists in $X$. This completes the proof.
\end{proof}


\begin{lemma}
\label{lem:3p}
If the size of blocks in $\B$ is $3p$ then $X$ 
has a Hamilton path. 
\end{lemma}

\begin{proof}
Note that $|\B|=2$ and $X_\B=K_2$. Let us denote the two
blocks of $\B$ by $B$ and $B'$.
We first show that in the case when $p = 5$ we can assume a $(6,5)$-semiregular
automorphism $\rho$ of $X$, with $\rho \in G$, to exist. Suppose on the contrary that $X$ does not
contain a Hamilton path and that no such $\rho \in G$ exists. By Proposition~\ref{pro:jack} the valency
of $X$ is at most $9$. Let $P \leq G$ be a 
Sylow $5$-subgroup of $G$. In view of Proposition~\ref{pro:wielandt} and 
Lemma~\ref{lem:blocksemi} $P$ has six orbits
of length $5$ on $X$. Denote them by $S_i$, $i \in \ZZ_6$. 
With no loss of generality assume that $S_i \subset B$ for $i = 0,1,2$. 
Proposition~\ref{pro:semireg}
implies that there exists some $\psi \in G$, such that $\psi|_B$ is $(3,5)$-semiregular. With no loss of generality
assume that $\psi$ is of order $5$ and $\psi\in P$. The orbits of $\psi$ on $B$ are then $S_0$, $S_1$ and $S_2$.
In view of our assumptions $\psi|_{B'}$ is not semiregular. Moreover, 
$\psi|_{B'} \neq Id$, as otherwise $\psi \alpha^{-1}\psi \alpha$ is $(6,5)$-semiregular on $X$, where
$\alpha \in G$ is such that $B\alpha = B'$.
Thus $\psi$ has at least one orbit of length $5$ on $B'$ and at least $5$ fixed points. We can assume that
this orbit of length $5$ is $S_3$ and that the $5$ fixed points are the vertices of $S_4$. 
As $X_\B \cong K_2$, we can assume that $S_1 \sim S_4$. Since $S_1$ and $S_4$ are orbits of $P$, 
it is clear that $[S_1,S_4] = K_{5,5}$. Moreover, since $x$ has
at most $9$ neighbors, the valency of $[B,B']$ is $5$, and so $[B,B'] = 3K_{5,5}$. 
Since $S_1$ is a subset of the block $B$,
it is now clear that $S_1$ itself is a block for $G$. Lemma~\ref{lem:blocksemi} implies that
the block system arising from $S_1$ coincides with $\{S_i \ |\ i \in \ZZ_6\}$. 
Using the fact that $X$ is connected one can see that there exist adjacent vertices $u$ and $v$ of $B'$
such that $\psi$ fixes precisely one of them. But then the valency of $X$ is at least $10$, a contradiction
which proves our claim.

Therefore, Proposition~\ref{pro:semireg} implies that 
we can assume that a $(6,p)$-semiregular automorphism
$\rho$ of $X$ such that $\rho \in G$ exists. Let $\S = \{S_i\ |\ i \in \ZZ_6\}$ be
the set of its orbits. By Lemma~\ref{lem:blocksemi} each block in $\B$ is a union of three orbits of $\rho$.
With no loss of generality we can assume that $B = S_0 \cup S_1 \cup S_2$ and $B' = S_3 \cup S_4 \cup S_5$.
In view of regularity of the bipartite graph $[B,B']$,  the subgraph $\bar{X}_\S$
of $X_\S$, which is obtained from $X_\S$ by deleting  the edges between the orbits inside the blocks $B$ and $B'$
(if they exist),
is clearly one of the graphs $Y_i$, $i \in \{0,1,2,3,4,5\}$ of Figure~\ref{fig:rho3p}.
However, for each of the graphs $Y_i$, $i \geq 2$, the following holds: if there exists a multiedge of $X_\rho$, 
then there exists a Hamilton cycle of $X_\S$ which contains an edge
corresponding to a multiedge of $X_\rho$. By Lemma~\ref{lem:cyclelift} we can thus assume that
no multiedge exists in $X_\rho$ except possibly when $\bar{X}_\S = Y_0$ or $\bar{X}_\S = Y_1$. 
Regularity of $X$ then implies that $Y_4$ and $Y_5$ are not possible.

\begin{figure}[htbp]
\begin{footnotesize}
\begin{center}
\includegraphics[width=1.00\hsize]{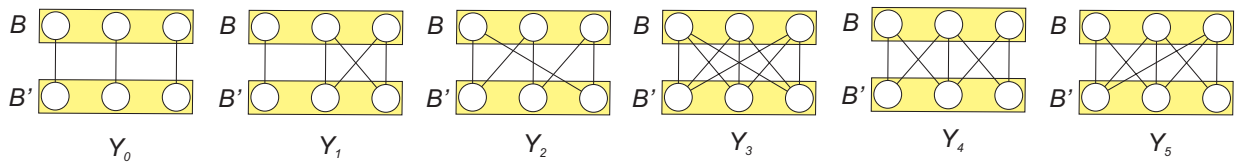}
\caption{\label{fig:rho3p} All possibilities for the subgraph  $\bar{X}_\S$ of $X_\S$
when $G$ has a complete block system $\B$ of two blocks of size $3p$.}
\end{center}
\end{footnotesize}
\end{figure}

If $X\la B \ra$ is a connected graph, then it contains a Hamilton cycle (as it is a vertex-transitive graph of 
order $3p$) and so $X$  has a Hamilton path in this case. 
We can thus assume that $X\la B \ra$ is not connected, and so 
it is isomorphic to $3pK_1$, to $pK_3$ or it is a disjoint union of three isomorphic connected $p$-circulants.
We consider each of the three cases separately. The technical details of each of them are left to the reader.

\setcounter{case}{0}
\begin{case} $X\la B  \ra \cong 3pK_1$. \end{case}
As $X$ is connected, the quotient graph $X_\S = \bar{X}_\S$ is one of
$Y_2$ and $Y_3$. If $X \cong Y_2\cong C_6$, then connectedness of $X$ and Lemma~\ref{lem:cyclelift} imply that
the lift of $Y_2$ contains a Hamilton cycle of $X$. 
If however $X_\S \cong Y_3\cong K_{3,3}$, then one can see that 
some Hamilton cycle of $X_\S$, whose lift contains a Hamilton cycle of $X$, exists.

\begin{case} $X\la B \ra \cong pK_3$. \end{case} Then of course also  $X\la B' \ra \cong pK_3$.
It is clear that each $K_3$ in $B$, $B'$ intersects all the orbits of $\rho $ in $B$ and $B'$, respectively.
Suppose first that $\bar{X}_\S \cong Y_0$. Then every edge of $X_\S$ is contained on some
Hamilton cycle of $X_\S$. Hence Lemma~\ref{lem:cyclelift} 
implies that we can assume that no multiedge exists in $X_\rho$. 
If there exists a Hamilton cycle of $X_\S$ whose lift contains a Hamilton cycle of $X$, we are done. If not,
the connectedness of $X$ implies that $X \cong C_{3p} \square K_2$, and so $X$ contains a Hamilton cycle.
If $\bar{X}_\S \cong Y_1$ then there exists a multiedge of $X_\rho$ 
that is contained in a Hamilton cycle of $X_\S$, and so a Hamilton cycle of $X$ exists.
Finally, if $\bar{X}_\S$ is isomorphic to $Y_2$ or to $Y_3$ it is easy to see that
some Hamilton cycle of $X_\S$, 
whose lift contains a Hamilton cycle of $X$, exists.

\begin{case} $X\la B \ra$ \end{case} is isomorphic to a disjoint union of three isomorphic connected
$p$-circulants. Then the quotient graph $X_\S = \bar{X}_\S$ is one of
$Y_2$ or $Y_3$, and so it has a Hamilton cycle. 
As the six $p$-circulants are precisely the graphs $X\la S_i\ra$, where $i \in \ZZ_6$, 
a Hamilton path exists in $X$. 
This completes the proof of Lemma~\ref{lem:3p}.
\end{proof}


\section{Quasiprimitive  graphs}
\label{sec:quasi}
\indent

Throughout this section let $X$ denote a connected quasiprimitive graph of order $6p$.
In \cite{MSZ95} a complete characterization of quasiprimitive graphs of order $pqr$, where $p$, $q$ and $r$ 
are distinct primes, was given via the well known generalized 
orbital graph constraction relative to certain simple groups 
having an imprimitive permutation representation of degree $pqr$. All the possible group
actions are given in Tables A and B in \cite[p.~298-299]{MSZ95}. For our purposes (we require that $pqr = 6p'$) 
only a handful of group actions needs to be considered. They
are given in Table~\ref{tab:qprim}.
Note that only row~11 of Table~\ref{tab:qprim} corresponds to an infinite family of actions giving rise
to quasiprimitive graphs of order $6p$. Lemma~\ref{lem:row5} shows that each of the quasiprimitive 
graphs corresponding
to an action from this infinite family has a Hamilton cycle.
As for the other rows of Table~\ref{tab:qprim}, each case is investigated separately.
More precisely, we consider all the possible generalized orbital graphs 
and study their structural properties (using program package {\sc Magma}~\cite{Mag})
which allows us to easily find a Hamilton path.
In fact, in all the graphs, except for the truncation of the Petersen graph, 
a Hamilton cycle is found. 

\begin{table}[htb] 
 $$
 \begin{array}{|c|c|c|}
   \hline
   \textrm{row}&p & \textrm{Action} \\\hline \hline
   1& 5 & A_5\; \textrm{on\; cosets\; of\; $\ZZ_2$} \\ \hline
   2& 7 & A_7\; \textrm{on\;  cosets\; of\; $A_5$}\\ \hline
   3&  11 & \PSL(2,11)\; \textrm{on\; cosets\; of\; $D_{10}$} \\ \hline
   4&  7 & \PSL(3,2)\; \textrm{on\; cosets\; of\; $\ZZ_4$} \\ \hline
   5&  7 & \PSL(3,2)\; \textrm{on\; cosets\; of\; $\ZZ_2^2$} \\ \hline
   6&  13 & \PSL(3,3)\; \textrm{on\; cosets\; of\; $\ZZ_3^2\rtimes \ZZ_8$} \\ \hline
   7&  13 & \PSL(3,3)\; \textrm{on\; cosets\; of\; $\ZZ_3^2\rtimes D_8$} \\ \hline
   8&  13 & \PSL(3,3)\; \textrm{on\; cosets\; of\; $\ZZ_3^2\rtimes Q_8$} \\ \hline
   9&  31 & \PSL(3,5)\; \textrm{on\; cosets\; of\; $\ZZ_5^2\rtimes (\ZZ_5\rtimes \ZZ_4^2)$} \\ \hline
   10&  5 & A_6\; \textrm{on\; cosets\; of\; $A_4$} \\ \hline
   11&  \frac{k+1}{2} & \PSL(2,k)\; \textrm{on\; cosets\; of\; $\ZZ_k\rtimes \ZZ_{(k-1)/6}$\; where $3\mid \frac{k-1}{2}$
   \; and\;  $k=s^m$} \\ \hline
   \end{array}
 $$
 \caption{\label{tab:qprim}  Actions giving rise to quasiprimitive graphs of order $6p$.}
 \end{table}

Let $G$ be a group acting on the cosets of its subgroup $H$ in a natural way. We say that the set $\O(G,H)$ of 
generalized orbital graphs (in short GOGs)
of this action is a {\em minimal connected orbital graph set} for this action
if each connected GOG corresponding to this action contains some graph of $\O(G,H)$ as a spanning subgraph.
As we are only interested in whether a given GOG contains a Hamilton path (or a Hamilton cycle)
Proposition~\ref{pro:jack} implies that we can disregard the graphs from $\O(G,H)$ whose valencies are
at least $[G:H]/3$. We let the remaining set of GOGs be the set $\R(G,H)$ of {\em relevant graphs} for
this action. It is now clear that in order to show that each GOG corresponding to the above mentioned action of $G$
contains a Hamilton path (Hamilton cycle) we only need to show that each GOG of $\R(G,H)$ has this property.

We now describe the method of obtaining $\R(G,H)$ for the action of row~1 of Table~\ref{tab:qprim} in full detail.
The other actions are dealt with in a similar way, so we only give the relevant graphs and leave
the details to the reader. Each relevant graph $X$ will be represented in a structural way given by some
semiregular automorphism $\varphi$ of $X$ from which the existence of a Hamilton cycle will be clear (except
for the truncation of the Petersen graph).
In the case when $\varphi$ is $(6,p)$-semiregular its symbol (for the definition see the next paragraph) 
will be given. In other cases we give the graph in its Frucht's notation.

Let $\rho$ be a $(6,p)$-semiregular automorphism and let
$S_i$, $i \in \ZZ_6$, be its orbits. 
Choose $s_i \in S_i$ and define the following subsets of $\ZZ_p$.
For $i,j \in \ZZ_6$, we let $R_{i,j} =\{r \in \ZZ_p\ |\ s_i \sim s_j\rho^r\}$. Note that $R_{j,i} =-R_{i,j}$.
It is clear that the collection of all $R_{i,j}$ completely determines $X$.
The $6\times 6$-matrix $\M_\rho(X) =(R_{i,j})_{i,j}$, whose $(i,j)$-th entry is
the set $R_{i,j}$, is the {\em symbol} of $X$ relative to ($\rho$, $s_0$, $s_1$, $s_2$, $s_3$, $s_4$, $s_5$).
 \medskip

\noindent
{\em Graphs corresponding to row~1 of Table~\ref{tab:qprim}}:
Note that these graphs are of order $30$. In the action of $A_5$ on the cosets of $\ZZ_2$  
we get that $\ZZ_2$  has $15$ nontrivial suborbits, 
$7$ of which are self-paired. Of the seven self-paired suborbits, 
six are of length $2$ and one is of length $1$. The non-self-paired suborbits are of length $2$.
Denote the $15$ nontrivial suborbits by $U_i$, $i \in \{1,2,\ldots , 15\}$, where $U_1$ is of length $1$,
$U_2, U_3, \ldots , U_7$ are the self-paired suborbits of length $2$ and $U_{2i}$ is paired with $U_{2i+1}$ for
$i \in \{4,5,6,7\}$.

\begin{figure}[htbp]
\begin{footnotesize}
\begin{center}
\includegraphics[width=0.70\hsize]{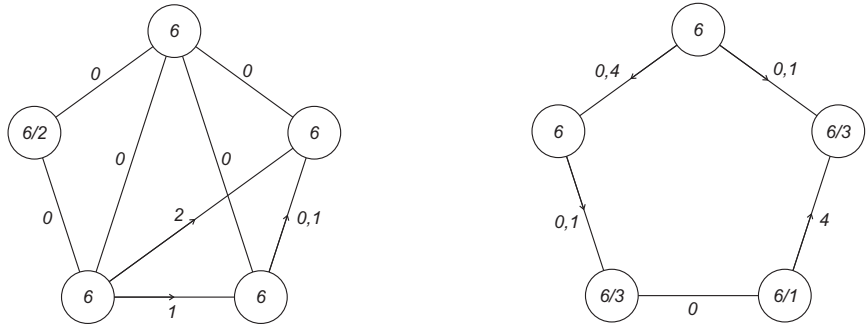}
\caption{\label{fig:qprim_row1} Two graphs given in the Frucht's notation under a $(5,6)$-semiregular 
automorphism.}
\end{center}
\end{footnotesize}
\end{figure}

The unions $U_{2i} \cup U_{2i+1}$, where $i \in \{4,5,6,7\}$, give rise to three nonisomorphic graphs, 
one of which is disconnected (with no loss of generality assume that this graph corresponds to $U_{14} \cup U_{15}$).
The other two are given in Frucht's notation under a $(5,6)$-semiregular automorphism in Figure~\ref{fig:qprim_row1}. 
Using an argument similar to the one in the proof of 
Lemma~\ref{lem:cyclelift} one can see that these two graphs both contain a Hamilton cycle. 

It turns out that the graph arising from $U_1 \cup U_{14} \cup U_{15}$ is still disconnected. 
The graphs arising from $U_i \cup U_{14} \cup U_{15}$, where $i \in \{2,3, \ldots , 7\}$, 
are all connected and isomorphic either to $X_1$ or to $X_2$ of Table~\ref{tab:qprim_row1}, and so 
Lemma~\ref{lem:cyclelift} implies that a Hamilton cycle exists in $X$.
Therefore, we now only have to consider the GOGs arising from unions of some suborbits 
from $\{U_1, U_2, \ldots , U_7\}$.
 
\begin{figure}[htbp]
\begin{minipage}[b]{0.48\linewidth}
\centering
\includegraphics[width=0.6\linewidth]{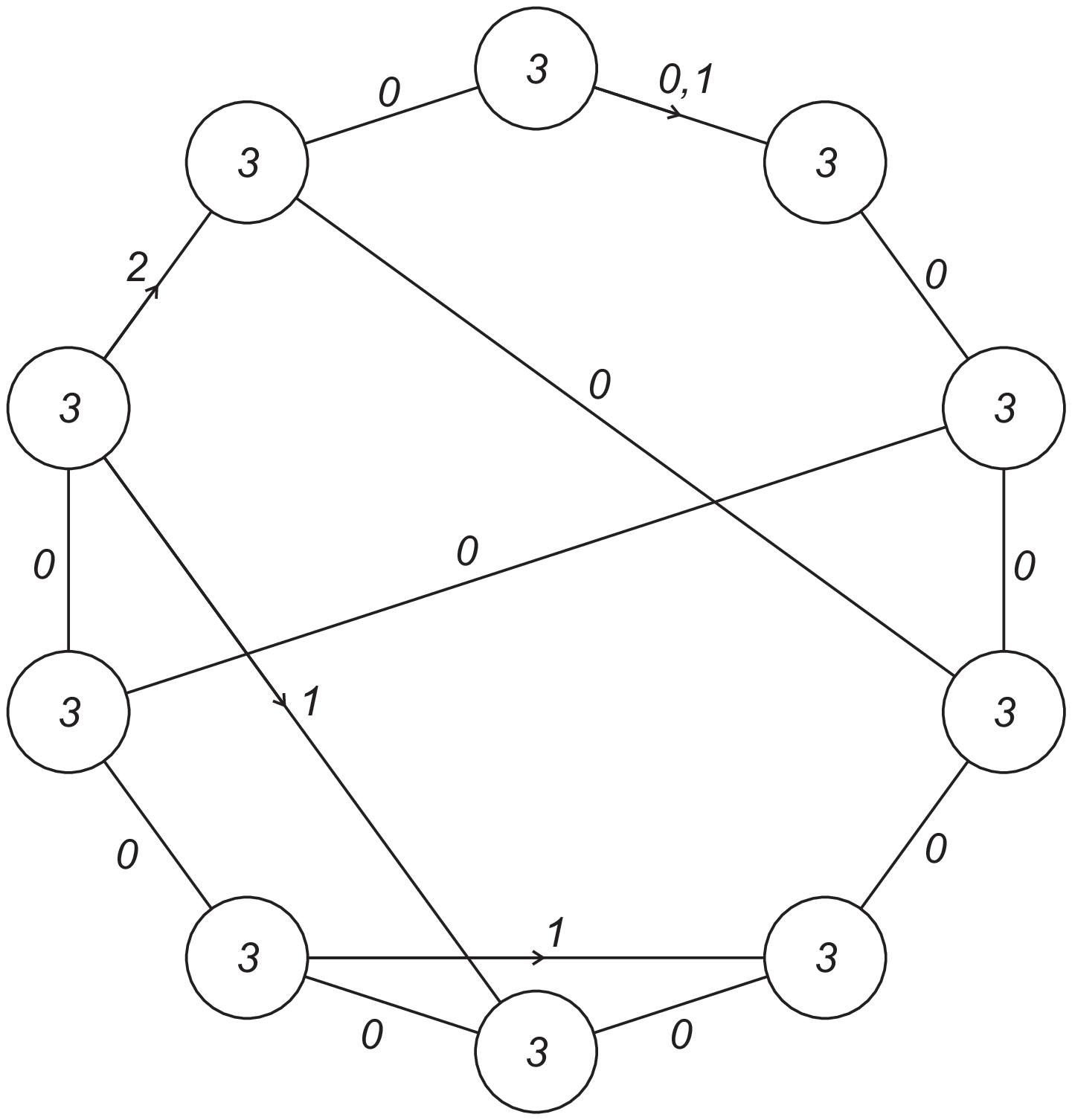}
\caption{\label{fig:qprim_row1a} A graph
 given in the Frucht's notation under a $(10,3)$-semiregular 
automorphism.}
\end{minipage}\hspace{0.46cm}
\begin{minipage}[b]{0.48\linewidth}
\centering
\includegraphics[width=0.6\linewidth]{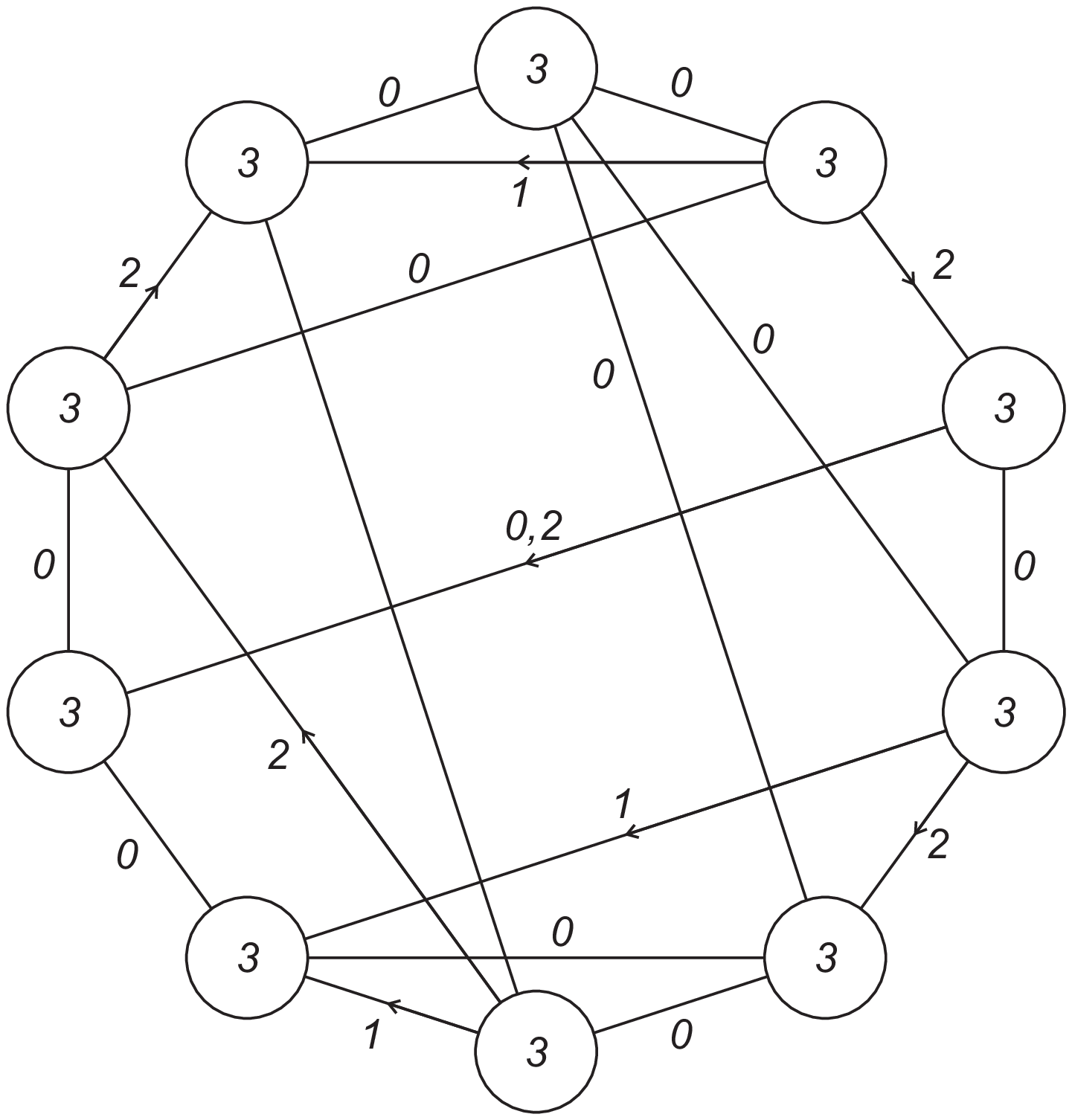}
\caption{\label{fig:qprim_row1b}  A graph
 given in the Frucht's notation under a $(10,3)$-semiregular 
automorphism. } 
\end{minipage}
\end{figure}

For every $i \in \{1,2,3, \ldots , 7\}$ the graph arising from the suborbit $U_i$ is disconnected, whereas
the graph arising from $U_1 \cup U_i$, $i \in \{2,3, \ldots , 7\}$, is connected and 
isomorphic either to the truncation of the Petersen graph, or to the graph of
Figure~\ref{fig:qprim_row1a} given in the Frucht's notation under a $(10,3)$-semiregular automorphism.
Lemma~\ref{lem:cyclelift} implies that the latter graph 
contains a Hamilton cycle.

Finally, the unions $U_i \cup U_j$, where $i,j \in \{2, 3, \ldots , 7\}$, 
give rise to five nonisomorphic connected graphs. These are the
graphs $X_3$, $X_4$, $X_5$ and $X_6$  of Table~\ref{tab:qprim_row1} and the graph of Figure~\ref{fig:qprim_row1b}
given in Frucht's notation under a $(10,3)$-semiregular automorphism.
Lemma~\ref{lem:cyclelift} implies that in all these cases the graph in question has a Hamilton cycle.

We have now clearly considered all the relevant graphs $\R(A_5, \ZZ_2)$.
Note also, that each GOG corresponding to the action of $A_5$ on the cosets of $\ZZ_2$ 
which contains the truncation of the Petersen graph as a proper 
spanning subgraph contains a Hamilton cycle.
We can thus conclude that each connected GOG arising from the action of $A_5$ on the cosets of $\ZZ_2$, except for 
the truncation of the Peterson graph, contains a Hamilton cycle.

\medskip

\noindent
{\em Graphs corresponding to row~2 of Table~\ref{tab:qprim}}: The relevant graphs are
given in Table~\ref{tab:qprim_row2}, and so it is clear that each GOG arising from this action contains a 
Hamilton cycle.

\medskip

\noindent
{\em Graphs corresponding to row~3 of Table~\ref{tab:qprim}}: The relevant graphs are
given in Table~\ref{tab:qprim_row3}, and so it is clear that each GOG arising from this action contains a 
Hamilton cycle.

\medskip

\begin{figure}[htbp]
\begin{minipage}[b]{0.48\linewidth}
\centering
\includegraphics[width=0.70\hsize]{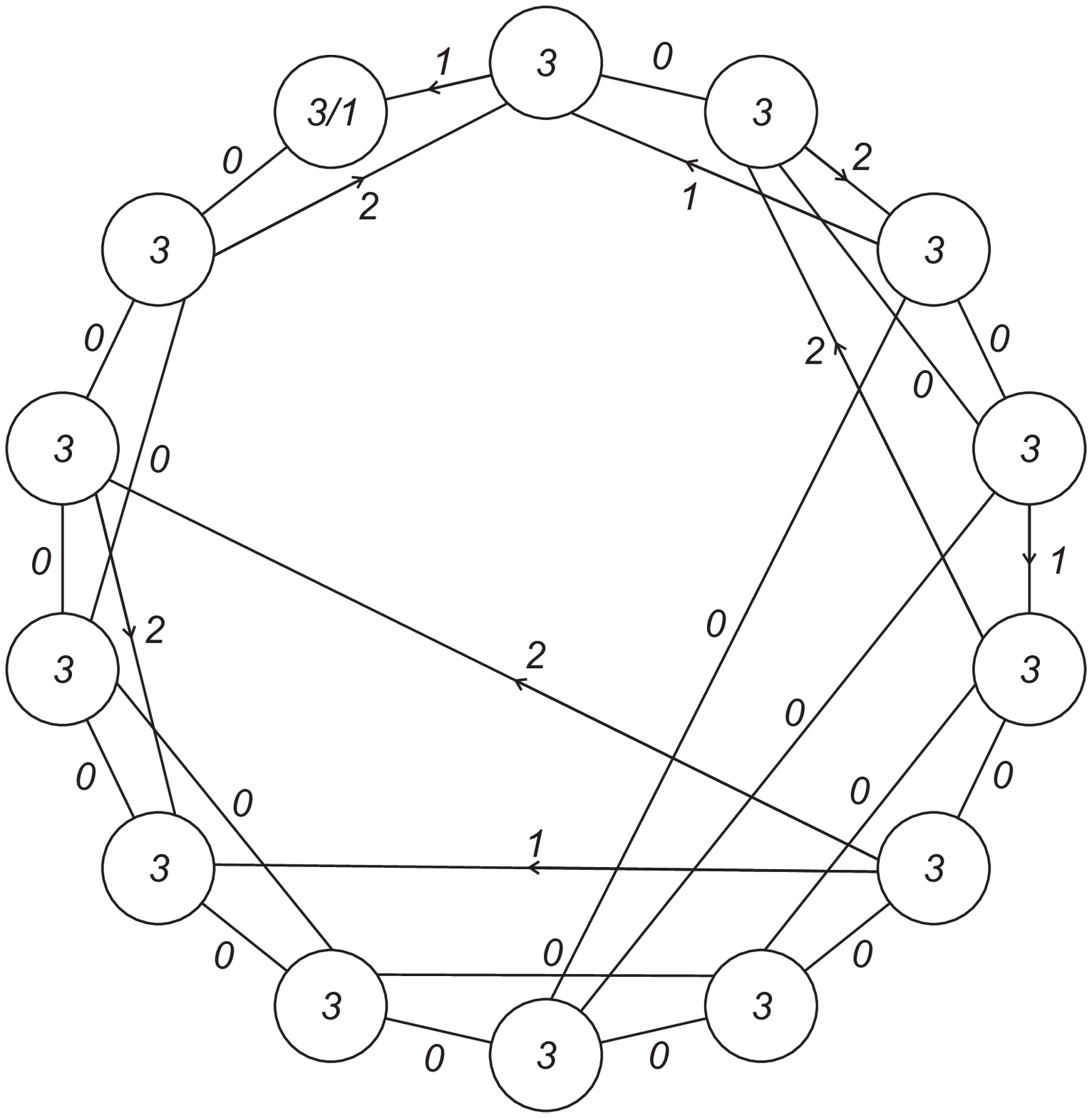}
\caption{\label{fig:qprim_row4} A graph 
given in the Frucht's notation under a $(14,3)$-semiregular 
automorphism.}
\end{minipage}\hspace{0.46cm}
\begin{minipage}[b]{0.48\linewidth}
\centering
\includegraphics[width=0.6\hsize]{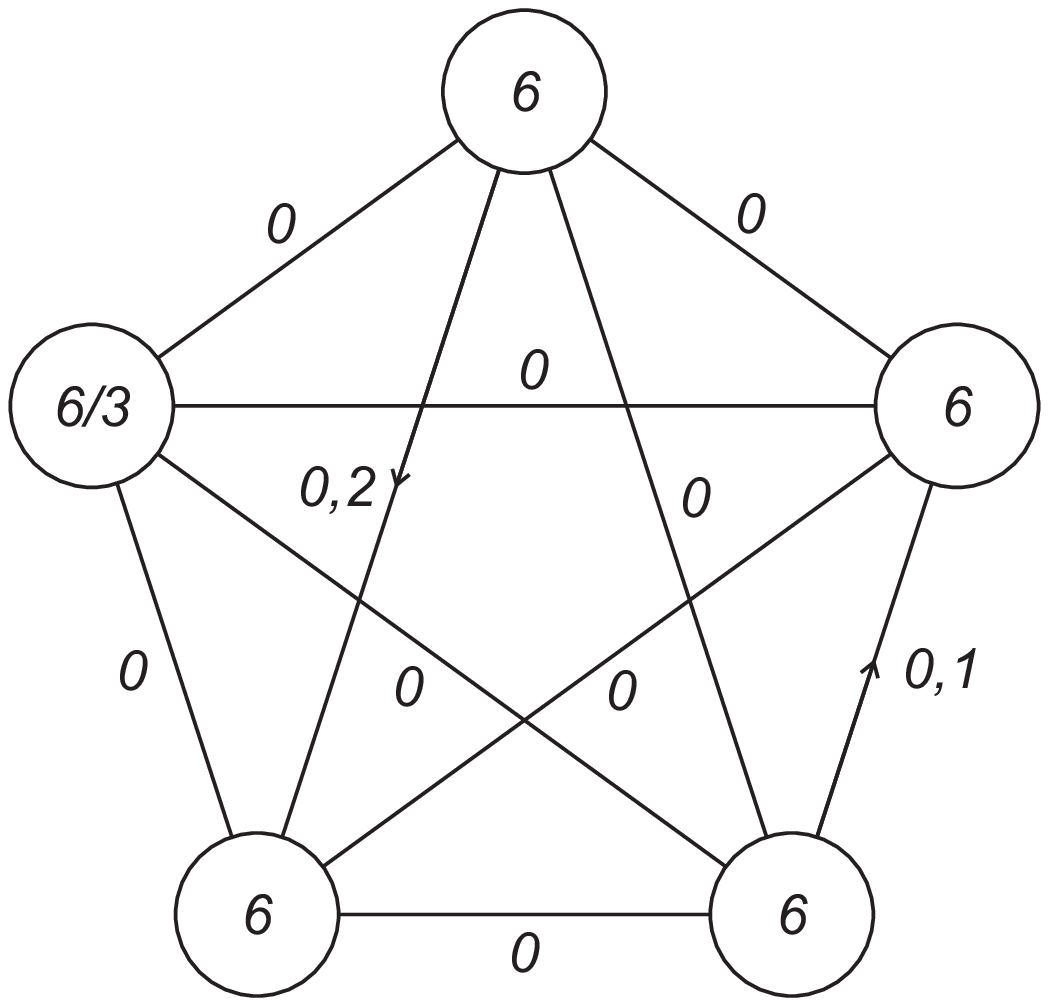}
\caption{\label{fig:qprim_row10} A graph
given in the Frucht's notation under a $(5,6)$-semiregular 
automorphism.}
\end{minipage}
\end{figure}

\noindent
{\em Graphs corresponding to row~4 of Table~\ref{tab:qprim}}: The relevant graphs are
given in Table~\ref{tab:qprim_row4} and Figure~\ref{fig:qprim_row4}, and so it is clear that each 
GOG arising from this action contains a 
Hamilton cycle. \medskip

\noindent
{\em Graphs corresponding to row~5 of Table~\ref{tab:qprim}}: The relevant graphs are
given in Table~\ref{tab:qprim_row5}, and so it is clear that each GOG arising from this action contains a 
Hamilton cycle. 

\medskip

\noindent
{\em Graphs corresponding to row~6 of Table~\ref{tab:qprim}}: It turns out that $\R(G,H) = \emptyset$ in this case,
and so each GOG arising from this action contains a Hamilton cycle.

\medskip

\noindent
{\em Graphs corresponding to row~7 of Table~\ref{tab:qprim}}: The relevant graphs are
given in Table~\ref{tab:qprim_row7}, and so it is clear that each GOG arising from this action contains a 
Hamilton cycle.

\medskip

\noindent
{\em Graphs corresponding to row~8 of Table~\ref{tab:qprim}}: It turns out that $\R(G,H) = \emptyset$ in this case,
and so each GOG arising from this action contains a Hamilton cycle.

\medskip

\noindent
{\em Graphs corresponding to row~9 of Table~\ref{tab:qprim}}: The relevant graphs are
given in Table~\ref{tab:qprim_row9}, and so it is clear that each GOG arising from this action contains a 
Hamilton cycle.

\medskip

\noindent
{\em Graphs corresponding to row~10 of Table~\ref{tab:qprim}}: The relevant graphs are
given in Table~\ref{tab:qprim_row10} and Figure~\ref{fig:qprim_row10}, and so 
 it is clear that each GOG arising from this action contains a 
Hamilton cycle.

\medskip

\noindent
{\em Graphs corresponding to row~11 of Table~\ref{tab:qprim}}: 
Lemma~\ref{lem:row5} below implies that each of the corresponding graphs contains a Hamilton cycle.

 \setcounter{case}{0}
 \begin{lemma}
 \label{lem:row5}
 Let $X$ be a graph corresponding to the action of row~11   of Table~\ref{tab:qprim}. Then 
 $X$ contains a Hamilton cycle.
 \end{lemma}
 
 \begin{proof}
 From \cite[Table~B and Section~4]{MSZ95} we can extract that the action of
 $G=\PSL(2,k)$ on the cosets of $\ZZ_k\rtimes \ZZ_{(k-1)/6}$ (the action of row~11
 of Table~\ref{tab:qprim}) gives rise to a vertex-transitive graph $X$ on which $G$
 has a complete block system $\B$ of $k+1 = 2p$ blocks of size $3$  with block stabilizer
 $G_B\cong \ZZ_k\rtimes \ZZ_{(k-1)/2}$.  Moreover, the permutation group $\bar{G}$ corresponding to
 the natural action of $G$ on $X_\B$ is doubly transitive, and so $X_\B$ is 
 isomorphic to the complete graph $K_{2p}$ and the bipartite graphs
 $[B,B']$, where $B,B'\in\B$, are all isomorphic. Note also, that $p \geq 7$.
  
 Since $3$ divides $\frac{k-1}{2}$ and $p = \frac{k+1}{2}$ it is clear that $p \equiv 1 \pmod 3$. Let $P$
 be a Sylow $3$-subgroup of $G = PSL(2,k)$ 
 and let $\bar{P}$ denote the
 permutation group corresponding to the natural action of $P$ on $X_\B$. 
 Since $\bar{P}$ is a $3$-group and
 $2p \equiv 2 \pmod 3$, there exist $B_0, B_1 \in \B$ which are fixed by $\bar{P}$. By Proposition~\ref{pro:wielandt}, however,
 $P$ acts transitively on each of the two blocks $B_0$ and $B_1$. Since $X_\B$ is a complete graph, there
 exist adjacent vertices $u \in B_0$ and $v \in B_1$. Let $\varphi \in P$ be an automorphism which does not fix $u$.
 If it fixes $v$, then $[B_0,B_1]$ is a complete bipartite graph $K_{3,3}$, and so $X$ is of valency at least 
 $6p - 3$, in which case Proposition~\ref{pro:jack} applies. 
 We can therefore assume that $\varphi$ does not  fix $v$. Then
 $[B_0,B_1]$ contains $3K_2$ as a subgraph. If $[B_0,B_1]$ is not isomorphic to $3K_2$ or if $X\la B_0 \ra$ is not an
 independent set, then the valency of $X$ exceeds $2p$, and we can again apply Proposition~\ref{pro:jack}. 
 
 We can now assume that $[B,B'] \cong 3K_2$ and $X\la B \ra = 3K_1$ for all $B,B' \in \B$. 
 As $p \geq 7$, Proposition~\ref{pro:semireg} implies that a $(6,p)$-semiregular automorphism $\rho$ of $X$, where 
 $\rho \in G$, exists. Denote its orbits by $\S = \{S_i\ |\ i \in \ZZ_6\}$.
 By Lemma~\ref{lem:blocksemi}, we have $|S_i \cap B| \in \{0,1\}$ for all $i \in \ZZ_6$ and $B \in \B$.
 It is clear that  we can then assume that $\mathcal{A} = S_0 \cup S_1 \cup S_2$ is a
 union of $p$ blocks from $\B$ and that $\mathcal{A}' = 
 S_3 \cup S_4 \cup S_5$ is a union of the other $p$ blocks from $\B$.
 In view of our assumptions each vertex in $\mathcal{A}$ has $p$ neighbors in $\mathcal{A}'$ and vice versa.
 Suppose there exists an orbit $S_i$, with no loss of generality assume it is $S_0$, such that
 $X\la S_0 \ra = K_p$ and $S_0$ is adjacent to only one of the orbits from $\mathcal{A}'$, say to $S_3$.
 This implies that $[S_0, S_3] = K_{p,p}$.
 Note that the vertices of $S_0$
 are characterized by the fact that they are adjacent to all the vertices of $S_0 \cup S_3$ (except to itself).
 Moreover, as $X$ is connected, each vertex of $S_3$ has at least one neighbor outside $S_0 \cup S_3$.
 It is now clear that $S_0$ is a block of imprimitivity for $G$.
 But this implies that the quotient graph corresponding to the imprimitivity block system arrising from
 $S_0$ is a vertex-transitive graph of order $6$ which thus contains a Hamilton cycle. It is now clear
 that $X$ also has a Hamilton cycle.
 
 We can thus assume that for each $S_i \in \S$ the following holds:
 if $X\la S_i \ra = K_p$ then the valency of $S_i$ in the subgraph $Y = [\mathcal{A}, \mathcal{A}']$ of
 $X_\S$ is at least two. 
 Note also that if $S_j$ is the only neighbor of $S_i$ in $Y$, then
 $[S_i,S_j] = K_{p,p}$, and so $S_i$ is the only neighbor of $S_j$ in $Y$ as well.
 We distinguish two cases depending on whether the graph $Y$ contains a vertex of valency $1$ or not.
 
 \begin{case}
 \end{case}
 There exists a vertex of $Y$ of valency one. 
 With no loss of generality assume that the only neighbor  of $S_0$ in $\mathcal{A}'$ is $S_3$.
 We distinguish two cases depending on the valency $d$ of $S_1$ in $Y$.
 
 If $d = 1$, say $S_1 \sim S_4$, then the valency of $S_2$ in $Y$ is also $1$,
 and so $S_2 \sim S_5$.
 In view of the above remarks each of $S_i \in \mathcal{A}$ has at least one neighbor inside $\mathcal{A}$, and
 the same holds for $\mathcal{A}'$.
 We can thus assume that $S_0 \sim S_1$, $S_0 \sim S_2$ and $S_3 \sim S_4$. Moreover, $S_5$ is adjacent to one of 
 $S_3$ and $S_4$. If $S_5 \sim S_3$, then $S_0S_1S_4S_3S_5S_2S_0$ is a Hamilton cycle of $X_\S$ which
 contains an edge corresponding to a multiedge of $X_\rho$, so Lemma~\ref{lem:cyclelift} applies.
 Suppose then that $S_5 \not\sim S_3$, and so $S_5 \sim S_4$. Note that this also implies that
 $X\la S_3 \ra \neq pK_1$ (otherwise the valency of the vertices of $S_4$ exceeds $2p-1$).
 It is clear that then a Hamilton path of $[S_0,S_3]$ with endvertices
 in $S_0$ exists. As $[S_1,S_4] \cong [S_2,S_5] \cong K_{p,p}$, $S_0 \sim S_1$ and $S_0 \sim S_2$, 
 the existence of a Hamilton cycle of $X$ is evident.
 
 If $d > 1$, then clearly $[S_1 \cup S_2, S_4 \cup S_5] \cong K_{2,2}$. As the valency of $S_0$ in $Y$ is
 one, we have $X\la S_0 \ra \not\cong K_p$, and so $S_0$ is adjacent to at least one of $S_1, S_2$.
 Similarly, $S_3$ is adjacent to at least one of $S_4, S_5$. It is easy to see that a Hamilton cycle of $X_\S$ 
 containing the edge $S_0S_3$ exists in this case, so Lemma~\ref{lem:cyclelift} applies.
 
 \begin{case}
 \end{case} No vertex of valency $1$ exists in  $Y$. It is straightforward to check that in this case
 a Hamilton cycle of $X_\S$ containing an edge corresponding to a multiedge of $X_\rho$ exists, and
 so Lemma~\ref{lem:cyclelift} applies. We leave the details to the reader. 
\end{proof}

 \begin{table}[h!]
{\tiny
\begin{minipage}[b]{0.48\linewidth}
\centering
 $$
 \begin{array}{|c||c|c|c|c|c|c|}  
\hline

& X_1&X_2&X_3&X_4&X_5&X_6\\ 
\hline
\hline
  p&5&5&5&5&5&5\\  
\hline
  |V(X_i)|& 30&30&30&30&30&30\\  
\hline
  val& 6&6&4&4&4&4\\  
\hline\hline
  R_{0,0}&\emptyset&\emptyset&\emptyset&\emptyset&\emptyset&\emptyset\\ 
\hline
  R_{1,1}& \emptyset&\emptyset&\emptyset&\pm 1&\pm 2&\pm 2\\ 
\hline
  R_{2,2}& \emptyset&\emptyset&\emptyset&\emptyset&\emptyset&\emptyset\\ 
\hline
  R_{3,3}& \pm 1&\emptyset&\emptyset&\emptyset&\emptyset&\emptyset\\ 
\hline
  R_{4,4}& \emptyset&\emptyset&\emptyset&\emptyset&\emptyset&\pm 1\\ 
\hline
  R_{5,5}& \emptyset&\emptyset&\emptyset&\emptyset&\emptyset&\emptyset\\  
\hline
  R_{0,1}& 0&0,1& 0&0& 0& 0\\ 
\hline
  R_{1,2}& 3&4& 0&0& 0& 2\\ 
\hline
  R_{2,3}& 4&0,3& 0,2&0& 0& 2\\ 
\hline
  R_{3,4}& 0&2& 2&2,4& 0,1& 2\\  
\hline
  R_{4,5}& 0,2&3& 4&4& 0,1& 0\\  
\hline
  R_{5,0}& 0,2&0& 0,1&0& 0& 0\\ 
\hline
  R_{0,2}& 0,1&0& \emptyset&0& 0,3& 0\\ 
\hline
  R_{1,3}& 0&\emptyset& 0&\emptyset& \emptyset& \emptyset\\  
\hline
  R_{2,4}& 4&2& 4&\emptyset& \emptyset& \emptyset\\ 
\hline
  R_{3,5}& \emptyset&0,3& \emptyset&3& 0&1\\ 
\hline
  R_{4,0}& \emptyset&0& 0&0& \emptyset& \emptyset\\ 
\hline
  R_{5,1}& 0&1& 0&\emptyset& \emptyset& \emptyset\\ 
\hline
  R_{0,3}& 0&0& \emptyset&\emptyset& \emptyset& 0\\ 
\hline
  R_{1,4}& 0,4&0,1& \emptyset&\emptyset& \emptyset& \emptyset\\ 
\hline
  R_{2,5}& 4&1& \emptyset&2& \emptyset& 2\\ 
\hline
\hline
 \end{array}
  $$
\caption{\label{tab:qprim_row1}   Relevant graphs  
  corresponding to the action of row~1 of Table~\ref{tab:qprim}.} 
\end{minipage}\hspace{0.3cm}
\begin{minipage}[b]{0.48\linewidth}
\centering
$$
 \begin{array}{|c||c|c|c|c|}  
\hline
 & X_1&X_2&X_3&X_4\\ 
\hline
\hline
  p&7&7&7&7\\  
\hline
 |V(X_i)|& 42&42&42&42\\   
\hline
  val& 10&6&10&6 \\  
\hline\hline
  R_{0,0}&\pm 1 & \emptyset& \emptyset& \emptyset\\ %
\hline
  R_{1,1}& \pm 3 & \emptyset& \emptyset&\pm 3\\ %
\hline
  R_{2,2}& \pm 2& \emptyset& \emptyset& \emptyset\\ %
\hline
  R_{3,3}& \pm 3& \emptyset& \emptyset&\pm 1\\ %
\hline
  R_{4,4}& \pm 1& \emptyset& \emptyset& \emptyset\\ %
\hline
  R_{5,5}& \pm 2& \emptyset& \emptyset& \pm 2\\ %
\hline
  R_{0,1}& 1,5 & 0,5& 0,1& 0\\ %
\hline
  R_{1,2}&4,5 & 0& 2,6& 1,6\\ %
\hline
  R_{2,3}&1,3 & 0& 0,6& 0,3\\ %
\hline
  R_{3,4}&2,5 & 0,3& 4,6& 2\\ %
\hline
  R_{4,5}&0,6 & 0& 0,6& 4\\ %
\hline
  R_{5,0}&0,3 & 0& 0,4& 0\\ %
\hline
  R_{0,2}& 0,1& 0& 0,2& \emptyset\\ %
\hline
  R_{1,3}& \emptyset& 0& 1,6& \emptyset\\ %
\hline
  R_{2,4}& 0,4& 0& 4,5& \emptyset\\ %
\hline
  R_{3,5}& 0,6& 0& 3,6& \emptyset\\ %
\hline
  R_{4,0}& \emptyset& 0& 0,3& 0,1,3\\ %
\hline
  R_{5,1}& 0,2& 0& 0,5& \emptyset\\ %
\hline
  R_{0,3}& 0,5& 0& 0,1& 0\\ %
\hline
  R_{1,4}& 2,4& 0& 0,3& 5\\ %
\hline
  R_{2,5}& \emptyset& 0,6& 3,5& 4,5\\ %
\hline
\hline
 \end{array}
  $$
  \caption{\label{tab:qprim_row2}  Relevant graphs  
  corresponding to the action of row~2 of Table~\ref{tab:qprim}.}
\end{minipage}
}
\end{table}

 \begin{table}[h!]
{\tiny
\centering
 $$
 \begin{array}{|c||c|c|c|c|c|c|}
\hline

& X_1& X_2&X_3&X_4 &X_5 & X_6\\ 
\hline
\hline
  p&11&11&11&11&11&11\\  
\hline
 |V(X_i)|& 66& 66 &66&66&66&66\\  
\hline
  val& 5 & 10 & 10& 10& 10& 20\\  
\hline\hline
  R_{0,0}& \emptyset & \pm 4 & \pm 2 & \pm 3,\pm  5 & \emptyset& \pm 1 \\ 
\hline
  R_{1,1}& \emptyset &  \pm 1 & \pm 4 & \pm 1,\pm 5&  \emptyset&  \pm 2\\ 
\hline
  R_{2,2}& \emptyset & \pm 3  &  \pm 3& \pm 1,\pm 2 &  \emptyset& \pm 5 \\ 
\hline
  R_{3,3}& \emptyset & \emptyset & \emptyset & \pm 2, \pm 4 & \emptyset& \pm 1,\pm 2,\pm 3,\pm 4, \pm 5 \\ 
\hline
  R_{4,4}& \emptyset &  \pm 2 &  \pm 5&  \pm 3,\pm 4 &\emptyset&  \pm 3 \\ 
\hline
  R_{5,5}& \emptyset &  \pm 5 & \pm 1 & \emptyset  &\emptyset&  \pm 4 \\  
\hline
  R_{0,1}& 0,7 &  0,3,7& 0 &  0,5 &  0,8&  0, \pm 1, 2\\ 
\hline
  R_{1,2}& 0,8 & 0,9,10 & 2 &  4,5 & 3,9&  \pm 1, 4,6 \\ 
\hline
  R_{2,3}& 9 &  4,8&  \pm 3&  2,4 &  0,1&  5,10\\ 
\hline
  R_{3,4}& 5 & 0,10 &  8,9&  6,10 &  1,9&  1,4\\  
\hline
  R_{4,5}& 1,2 & 0,\pm 2 & 6 &  5,8& 0,8 & 0,3,7,10 \\  
\hline
  R_{5,0}& 0,9 &  0,5,10& 0 & 0,6 &  0,5&  0,1,7,8\\ 
\hline
  R_{0,2}& \emptyset & \emptyset& 0,4  & \emptyset & 0,9 & 0,1,\pm 5 \\ 
\hline
  R_{1,3}& 2 &  2,8& 6,9 &  \emptyset&  3,10&  0,9\\  
\hline
  R_{2,4}& 0,6 & 0,\pm 3& 3& \emptyset & \pm 2 & 0,3,6,9 \\ 
\hline
  R_{3,5}& 1 & 1,9& 8,10& 0,7 &  1,6&  0,4\\ 
\hline
  R_{4,0}& \emptyset & \emptyset & 0,1& 0,8 & 0,2 & 0,7,8,10\\ 
\hline
  R_{5,1}& \emptyset & \emptyset & \pm 1& \pm 5 & 0,2 & 0,\pm 2 ,7 \\ 
\hline
  R_{0,3}& 0 & 0,2& 0,4&  \emptyset& 0,10 & 0,10 \\ 
\hline
  R_{1,4}& \emptyset & \emptyset& 1,9&  \emptyset& 0,1 & \pm 1,2,4 \\ 
\hline
  R_{2,5}& \emptyset &  \emptyset& 1,6&  0,2&  2,6& 3,5,9,10 \\ 
\hline
\hline
 \end{array}
  $$
\caption{\label{tab:qprim_row3}   Relevant graphs  
  corresponding to the action of row~3 of Table~\ref{tab:qprim}.} 
}
\end{table}

 \begin{table}[h!] 
{\tiny
\begin{minipage}[b]{0.48\linewidth}
\centering
$$
 \begin{array}{|c||c|c|c|c|}  
\hline
 & X_1 & X_2 & X_3& X_4\\ 
\hline
\hline
  p&7 & 7& 7& 7\\  
\hline
  |V(X_i)|&42 & 42& 42& 42\\   
\hline
  val& 8 & 8 & 8 &   8 \\  
\hline\hline
  R_{0,0}& \pm 1 & \pm 2& \emptyset& \emptyset  \\ %
\hline
  R_{1,1}& \pm 1 & \pm 3& \emptyset&\emptyset  \\ %
\hline
  R_{2,2}& \pm 2 &\pm 1& \emptyset& \emptyset \\ %
\hline
  R_{3,3}& \pm 2 &\pm 1& \emptyset& \emptyset \\ %
\hline
  R_{4,4}& \pm 3 &\pm 3& \emptyset& \emptyset \\ %
\hline
  R_{5,5}& \pm 3 &\pm 2& \emptyset& \emptyset \\ %
\hline
  R_{0,1}& 0,4 & 0,3& 0& 0,5 \\ %
\hline
  R_{1,2}&  0,3&0,6& \pm 1 & 0,3  \\ %
\hline
  R_{2,3}&  5,6&4,5& 0,1 & 1,2 \\ %
\hline
  R_{3,4}&  4,5&0,1&  1 &  4,5\\ %
\hline
  R_{4,5}&  2,4&0,4&  0,5& 0,3 \\ %
\hline
  R_{5,0}&  0,2&0,5&  0,6& 0,5 \\ %
\hline
  R_{0,2}&  \emptyset&0,2& 0& 0,1  \\ %
\hline
  R_{1,3}&  \emptyset&\emptyset& 2,6& 2,4  \\ %
\hline
  R_{2,4}&  \emptyset&\emptyset& 5&  \emptyset \\ %
\hline
  R_{3,5}&  \emptyset&0,5& 0&  0,5 \\ %
\hline
  R_{4,0}& \emptyset&\emptyset& 0,4&  0,1 \\ %
\hline
  R_{5,1}& \emptyset&\emptyset& 3 &  \emptyset\\ %
\hline
  R_{0,3}& 0,4&\emptyset& 0,1 & \emptyset \\ %
\hline
  R_{1,4}& 0,2& 3,6 & 0,3 &  2,6\\ %
\hline
  R_{2,5}& 2,3& \emptyset& 3,5  &  2,6\\ %
\hline
\hline
 \end{array}
  $$
  \caption{\label{tab:qprim_row4} Relevant graphs  
  corresponding to the action of row~4 of Table~\ref{tab:qprim}.}
\end{minipage}\hspace{0.3cm}
\begin{minipage}[b]{0.48\linewidth}
\centering
 $$
 \begin{array}{|c||c| c|c|c|}
\hline

& X_1&X_2&X_3&X_4\\ 
\hline
\hline
  p&7&7&7&7\\  
\hline
  |V(X_i)|& 42 & 42& 42&42\\  
\hline
  val& 8 & 5&8&6\\  
\hline\hline
  R_{0,0}& \pm 3& \emptyset & \emptyset& \emptyset\\ 
\hline
  R_{1,1}& \pm 1& \emptyset& \emptyset& \emptyset \\ 
\hline
  R_{2,2}&  \pm 2& \emptyset& \emptyset& \emptyset\\ 
\hline
  R_{3,3}& \pm 3 & \emptyset& \emptyset& \emptyset\\ 
\hline
  R_{4,4}& \pm 1 & \emptyset& \emptyset& \emptyset\\ 
\hline
  R_{5,5}& \pm 2& \emptyset& \emptyset& \emptyset\\  
\hline
  R_{0,1}& 0 & 0& 0& 0\\ 
\hline
  R_{1,2}& 3,4 &0,4 & 0,4& 0,2\\ 
\hline
  R_{2,3}& 5 &6& 6& 3\\ 
\hline
  R_{3,4}&  5& 0,1& 2& 4,5\\  
\hline
  R_{4,5}& 4,5 &5& 0,1&0\\  
\hline
  R_{5,0}& 0 &0,2& 0&0\\ 
\hline
  R_{0,2}& 0,2 &0& 0,1&\emptyset\\ 
\hline
  R_{1,3}& 0,4&6& 2,3&3\\  
\hline
  R_{2,4}& 3& \emptyset& 4,6&0\\ 
\hline
  R_{3,5}& \pm 1&6& 0,3&\emptyset\\ 
\hline
  R_{4,0}& 0,3& 0& 0,3&0,2\\ 
\hline
  R_{5,1}& 0& \emptyset& 2,4&0,1\\ 
\hline
  R_{0,3}& \emptyset& \emptyset& 0,2& 0,3\\ 
\hline
  R_{1,4}& \emptyset& 0& 3& \emptyset\\ 
\hline
  R_{2,5}& \emptyset& 5& 5& 0,4\\ 
\hline
\hline
 \end{array}
  $$
\caption{\label{tab:qprim_row5}  Relevant graphs  
  corresponding to the action of row~5 of Table~\ref{tab:qprim}.} 
\end{minipage}

}
\end{table}

 \begin{table}[h!] 
{\tiny
\begin{minipage}[b]{0.25\linewidth}
\centering
$$
 \begin{array}{|c||c|}
\hline
 & X_1\\ 
\hline
\hline
  p&13\\  
\hline
  |V(X_i)|& 78\\   
\hline
  val&  18\\  
\hline\hline
  R_{0,0}& \pm 5 \\ %
\hline
  R_{1,1}& \pm 6 \\ %
\hline
  R_{2,2}& \pm 2 \\ %
\hline
  R_{3,3}& \pm 1\\ %
\hline
  R_{4,4}& \pm 3\\ %
\hline
  R_{5,5}& \pm 4\\ %
\hline
  R_{0,1}& 0,7,8\\ %
\hline
  R_{1,2}& 0,7,11\\ %
\hline
  R_{2,3}& 0,11,12\\ %
\hline
  R_{3,4}& 1,10,11\\ %
\hline
  R_{4,5}& 3,4,7\\ %
\hline
  R_{5,0}& 0,4,8\\ %
\hline
  R_{0,2}& 0,2,5\\ %
\hline
  R_{1,3}& 4,5,10,11\\ %
\hline
  R_{2,4}& 0,8,11\\ %
\hline
  R_{3,5}& 1,4,5\\ %
\hline
  R_{4,0}& 0,3,8,11\\ %
\hline
  R_{5,1}& 2,4,8\\ %
\hline
  R_{0,3}& 0,5,12\\ %
\hline
  R_{1,4}& 5,8,11\\ %
\hline
  R_{2,5}& 0,2,4,11\\ %
\hline
\hline
 \end{array}
  $$
  \caption{\label{tab:qprim_row7} Relevant graph  
  corresponding to the action of row~7 of Table~\ref{tab:qprim}.}
\end{minipage}\hspace{0.3cm}
\begin{minipage}[b]{0.35\linewidth}
\centering
 $$
 \begin{array}{|c||c|c|}
\hline

& X_1& X_2\\ 
\hline
\hline
  p&31 &31 \\  
\hline
  |V(X_i)|& 186& 186\\  
\hline
  val& 10 & 50\\ 
\hline\hline
  R_{0,0}&  \emptyset& \pm 1, \pm 2, \pm 6, \pm 9, \pm 13 \\ 
\hline
  R_{1,1}& \emptyset & \pm 1, \pm 3, \pm 5, \pm 10, \pm 14\\ 
\hline
  R_{2,2}& \emptyset &  \pm 4, \pm 11, \pm 12, \pm 13, \pm 14 \\ 
\hline
  R_{3,3}& \emptyset & \pm 4, \pm 7, \pm 9, \pm 10, \pm 15 \\ 
\hline
  R_{4,4}& \emptyset &  \pm 2, \pm 3, \pm 7, \pm 8 , \pm 11\\ 
\hline
  R_{5,5}&  \emptyset & \pm 5, \pm 6, \pm 8, \pm 12, \pm 15\\  
\hline
  R_{0,1}&  0,9 & 0,5,8,12,16,20,24,28 \\ 
\hline
  R_{1,2}& 0,4 & 5,12,18,20,21,27,28,30 \\ 
\hline
  R_{2,3}& 0,20& 3,5,6,8,11,13,14,16 \\ 
\hline
  R_{3,4}& 0,28&  1,2,\pm 3, 14,\pm 15, 20\\  
\hline
  R_{4,5}&  0,26&  5,6,\pm 8, \pm 9, 18,27\\  
\hline
  R_{5,0}& 0,6& 0,\pm 7, \pm 10, 11,14,28 \\ 
\hline
  R_{0,2}& 0,13&  0,3,8,14,19,22,24,29 \\ 
\hline
  R_{1,3}& 0,24& 8,\pm 9, 14,16,20,27,28 \\  
\hline
  R_{2,4}& 0,17&  1,4,10,\pm 11, 16,25,26\\ 
\hline
  R_{3,5}& 0,23&  \pm 4, \pm 13, \pm 14, \pm 15\\ 
\hline
  R_{4,0}& 0,1&  0,4,9,\pm 10, 14,25,26\\ 
\hline
  R_{5,1}& 0,15& 1,3,10,12,14,16,23,25  \\ 
\hline
  R_{0,3}& 0,2&  0,10,12,15,18,23,26,29 \\ 
\hline
  R_{1,4}& 0,21&  2,\pm 4, 8,14,18,20,26 \\ 
\hline
  R_{2,5}&  0,12& 4,5,6,8,11,13,14,15 \\ 
\hline
\hline
 \end{array}
  $$
\caption{\label{tab:qprim_row9}   Relevant graphs  
  corresponding to the action of  row~9 of Table~\ref{tab:qprim}.} 
\end{minipage}
\hspace{0.3cm}
\begin{minipage}[b]{0.3\linewidth}
\centering
 $$
 \begin{array}{|c||c|c|}
\hline

& X_1 & X_2 \\ 
\hline
\hline
  p&  5 &5 \\  
\hline
  |V(X_i)|& 30 & 5\\  
\hline
  val&  8& 8\\ 
\hline\hline
  R_{0,0}&  \pm 2 & \emptyset \\ 
\hline
  R_{1,1}& \emptyset &  \emptyset \\ 
\hline
  R_{2,2}&  \emptyset&  \pm 1, \pm 2 \\ 
\hline
  R_{3,3}& \pm 1 & \emptyset \\ 
\hline
  R_{4,4}& \pm 2& \emptyset\\ 
\hline
  R_{5,5}&  \pm 1 & \pm 1,\pm 2\\  
\hline
  R_{0,1}&  0 &  0,1\\ 
\hline
  R_{1,2}&  0,\pm 1, 2& 0 \\ 
\hline
  R_{2,3}& 4 & 0 \\ 
\hline
  R_{3,4}&  3,4 & 1,3\\  
\hline
  R_{4,5}&  0,2 & 0\\  
\hline
  R_{5,0}&  0,4 & 0\\ 
\hline
  R_{0,2}&  0 & 0 \\ 
\hline
  R_{1,3}&  3& 0,1 \\  
\hline
  R_{2,4}& 4  & 1 \\ 
\hline
  R_{3,5}&  \emptyset & 3  \\ 
\hline
  R_{4,0}&   \emptyset & 0,4 \\ 
\hline
  R_{5,1}&  1& 1 \\ 
\hline
  R_{0,3}&  0,2& 0,2 \\ 
\hline
  R_{1,4}&  0 & \pm 1\\ 
\hline
  R_{2,5}&  2 &  \emptyset\\ 
\hline
\hline
 \end{array}
  $$
\caption{\label{tab:qprim_row10}  Relevant graphs  
  corresponding to the action of row~10 of Table~\ref{tab:qprim}.} 
\end{minipage}
}
\end{table}

In view of the fact that the connected vertex-transitive graphs of orders $4p$ and $2p^2$ contain a 
Hamilton cycle (except for the Coxeter graph) (see \cite{KM??, DM87}), the results of this section imply
that the following proposition holds.

\begin{proposition}
\label{pro:quasi}
	A connected quasiprimitive graph of order $6p$, $p$ a prime, which is not
	isomorphic to the truncation of the Petersen graph, contains a Hamilton cycle.
\end{proposition}

\clearpage



\section{Primitive  graphs}
\label{sec:pri}
\indent

Throughout this section let $X$ denote a primitive graph of order $6p$.
In \cite{GP} the complete characterization of possible primitive graphs of order $2pq$, where $p$ and $q$ are distinct odd primes, 
was given. Extracting the information about graphs of order $6p$ we find that the only primitive
graphs of order $6p$, $p$ a prime, are the ones arising from the actions given in Table~\ref{tab:prim}.
Below we  show that each of the corresponding graphs has a Hamilton cycle.
We let the GOGs and the relevant graphs corresponding to some action be defined as in Section~\ref{sec:quasi}.

\begin{table}[h!] 
$$
\begin{array}{|c|c|c|}
  \hline
  \textrm{row}&p & \textrm{Action\; of\; $\Aut X$} \\\hline \hline
  1& 17 & \PSL(2,17)\; \textrm{on\; cosets\; of\; $S_4$} \\ \hline
  2& 11 & S_{12}\; \textrm{on\; pairs}\\ \hline
  3&  31 & \PSL(3,5)\; \textrm{on\; cosets\; of\; $P_{1,2}$} \\ \hline
\end{array}
$$
\caption{\label{tab:prim} Primes $p$ for which there exists a graph $X$ on $6p$ vertices such that $\Aut X$
and all vertex-transitive subgroups of $\Aut X$ act primitively on $X$.}
\end{table}

\medskip

\noindent
{\em Graphs corresponding to row~1 of Table~\ref{tab:prim}}: 
The relevant graphs are the so called $H$-graph (see \cite{TRG}), which 
by~\cite{TRG} has a Hamilton cycle,  and the graphs isomorphic to one of the graphs 
$X_1$, $X_2$, $X_3$, $X_4$, $X_5$ and $X_6$ of Table~\ref{tab:prim_row1}. 
It is therefore clear that each GOG arising from this action contains a 
Hamilton cycle.
\medskip

\noindent
{\em Graphs corresponding to row~2 of Table~\ref{tab:prim}}:
Note that $X$ is of order $66$. 
If $\{1,2\}\sim \{i,j\}$, where $\{1,2\}\cap\{i,j\}=\emptyset$ then the valency of $X$ is at least
$45$, so Proposition~\ref{pro:jack} applies. Therefore, the neighbors set of $\{1,2\}$ is the set
$\{\{i,j\}\mid i\in\{1,2\}, j \in \{3,4,\ldots,12\}\}$. It turns out that under
the $(6,11)$-semiregular automorphism $(1\,2\,3\,4\,5\,6\,7\,8\,9\,10\,11)\in S_{12}$ the symbol of $X = X_7$
is as in Table~\ref{tab:prim_row1}. Lemma~\ref{lem:cyclelift} implies that a Hamilton cycle exists in $X_7$.
\medskip

\noindent
{\em Graphs corresponding to row~3 of Table~\ref{tab:prim}}: 
The relevant graphs are isomophic to the graphs of Table~\ref{tab:qprim_row9}
and so it is clear that each GOG arising from this action contains a 
Hamilton cycle.

 \begin{table}[h!]
{\tiny
 $$
 \begin{array}{|c||c|c|c|c|c|c|c|} 
\hline

& X_1&X_2&X_3&X_4&X_5&X_6&X_7\\  
\hline
\hline
  p&17&17&17&17&17&17&11\\ 
\hline
  |V(X_i)|&102&102&102&102&102&102&66\\  
\hline
  val& 6&8&12&24&24&24&20\\  
\hline\hline
  R_{0,0}&\emptyset& \pm 1 &\pm 5 &\pm 1,\pm 4&\pm 6&\pm 5, \pm 8& \pm 1\\ 
\hline
  R_{1,1}&\pm 7& \pm 2, \pm 8 &\pm 1, \pm 4&\pm 8 &\pm 3, \pm 5& \pm 1, \pm 7& \pm 2\\ 
\hline
  R_{2,2}&\pm 6& \pm 4&\pm 7&\pm 2& \pm 6, \pm 7& \pm 2, \pm 3& \pm 3\\ 
\hline
  R_{3,3}&\emptyset& \pm 2&\pm 3&\pm 2,\pm 8& \pm 3& \pm 4, \pm 6& \pm 4\\ 
\hline
  R_{4,4}&\pm 3& \pm 1,\pm 4&\pm 2, \pm 8&\pm 1& \pm  7& \emptyset& \pm 5\\  
\hline
  R_{5,5}&\pm 5 & \pm 8&\pm 6&\pm 4& \pm 5& \emptyset& \pm 1, \pm 2, \pm 3, \pm 4,\pm 5 \\  
\hline
  R_{0,1}&0& 0,8&0,16 &0,5,\pm 8, 13,14& 0,2,10,12& 0,1,4,8,10,14 & 0,1,9,10\\  
\hline
  R_{1,2}&0& 0,2&0,1& 2,11,14,16& 1,\pm 2,\pm 3,5,11,13 & 0,3,\pm 8, 11,12& 0,2,8,10\\  
\hline
  R_{2,3}&4,10& 11,15&1& 0,3,8,12 & 10,11,15,16& 3,4,5,9,10,11& 0,3,7,10\\  
\hline
  R_{3,4}&0& 0,16&11,13& 3,4,\pm 5, 10,11 & 0,5,6,16& \pm 6, 9,14&0,4,6,10\\  
\hline
  R_{4,5}& 16 & 2,15 &0,8& 4,5,11,15& 0,6,9,14& \pm 2, 4,5,\pm 7, 8, 14 & 0,5\\  
\hline
  R_{5,0}&0,5& 0,16&0& 0,6,7,16& 0,3,7,13& 0,1,7,11 & 0,10\\  
\hline
  R_{0,2}&0& \emptyset &0& 0,\pm 1,2,4,14 & 0,11& \emptyset & 0,1,8,9\\  
\hline
  R_{1,3}&2,12& \emptyset&8,12& 0,2,14,16& 11,14& \emptyset& 0,2,7,9\\  
\hline
  R_{2,4}&\emptyset& \emptyset&12,14& 3,5,13,15& 2,12& \pm 4, \pm 8& 0,3,6,9\\  
\hline
  R_{3,5}&16& \emptyset&16& 2,5,9,11,13,15& 3,4,\pm 5,6,7,8,16& 0,7,8,15& 0,4\\  
\hline
  R_{4,0}&0,3& \emptyset &0,8& 0,4,11,15& 0,\pm 1,\pm 8,10,12,14& 0,13,14,16& 0,4,5,10\\ 
\hline
  R_{5,1}&\emptyset& \emptyset&6,10& \pm 4, 7,15& 5,10& \pm 3,5,16& 0,9\\  
\hline
  R_{0,3}&\emptyset& 0,2 &0,2,5,14& \emptyset& 0,3,5,15& 0,1,3,5,13,15& 0,1,7,8\\  
\hline
  R_{1,4}&\emptyset& \emptyset&\emptyset& 4,5,9,10& 5,8,14,16& \pm 6,9,14& 0,2,6,8\\  
\hline
  R_{2,5}&\emptyset& 0,8 &\pm 2, \pm 8& 0,10,13,14& 1,\pm 2,5& 4,8,14,15& 0,3\\  
\hline
\hline
 \end{array}
  $$
  \caption{\label{tab:prim_row1} {Relevant graphs  
  corresponding to the action of row~1 and row~2 of Table~\ref{tab:prim}.}}}
\end{table}


The results of this section imply that the following proposition holds.

\begin{proposition}
\label{pro:prim}
	A primitive graph of order $6p$, $p$ a prime, contains a Hamilton cycle.
\end{proposition}


\section{The proof of the main theorem}
\label{sec:proof}
\indent

\noindent
\proofT
In view of the results from \cite{KM??,DM87}, we can assume that $p \geq 5$.
If $X$ is not genuinely imprimitive, then either Proposition~\ref{pro:quasi} or Proposition~\ref{pro:prim}
applies. If $X$ is genuinely imprimitive, then apply one of Lemma~\ref{lem:2}, Lemma~\ref{lem:3}, 
Lemma~\ref{lem:p}, Lemma~\ref{lem:6}, Lemma~\ref{lem:2p} and Lemma~\ref{lem:3p}, depending on the
size of the corresponding blocks. \hfill \Qed

\end{document}